\documentclass{amsart}
\usepackage{longsalch, xy,tikz,float}
\xyoption{all}
\usepackage{MnSymbol}
\usepackage{cleveref}
\title{Structure and cohomology of moduli of formal modules.}
\author{A. Salch}

\begin{document}
\begin{abstract}
Given a commutative ring $A$, a ``formal $A$-module'' is a formal group equipped with an action of $A$. There exists a classifying ring $L^A$ of formal $A$-modules. This paper proves structural results about $L^A$ and about the moduli stack $\mathcal{M}_{fmA}$ of formal $A$-modules. We use these structural results to aid in explicit calculations of flat cohomology groups of $\mathcal{M}_{fmA}^{2-buds}$, the moduli stack of formal $A$-module $2$-buds. For example, we find that a generator of the group $H^1_{fl}(\mathcal{M}_{fm\mathbb{Z}}; \omega)$, which also generates (via the Adams-Novikov spectral sequence) the first stable homotopy group of spheres, also yields a generator of the $A$-module $H^1_{fl}(\mathcal{M}_{fmA}^{2-buds}; \omega)$ for any torsion-free Noetherian commutative ring $A$. We show that the order of the $A$-modules $H^1_{fl}(\mathcal{M}_{fmA}^{2-buds}; \omega)$ and $H^2_{fl}(\mathcal{M}_{fmA}^{2-buds}; \omega\otimes \omega)$ are each equal to $2^{N_1}$, where $N_1$ is the leading coefficient in the $2$-local zeta-function of $\Spec A$. We also find that the cohomology of $\mathcal{M}_{fmA}^{2-buds}$ is closely connected to the delta-invariant and syzygetic ideals studied in commutative algebra: $H^0_{fl}(\mathcal{M}_{fmA}^{2-buds}; \omega\otimes \omega)$ is the delta-invariant of the largest ideal of $A$ which is in the kernel of every ring homomorphism $A\rightarrow \mathbb{F}_2$, and consequently $H^0_{fl}(\mathcal{M}_{fmA}^{2-buds}; \omega\otimes \omega)$ vanishes if and only if $A$ is a ring in which that ideal is syzygetic.
\end{abstract}
\maketitle
\tableofcontents

\section{Introduction and review of some known facts.}

\subsection{Introduction.}

\subsubsection{Formal modules.}
This paper is about the classifying ring $L^A$ and classifying Hopf algebroid $(L^A,L^AB)$ of formal $A$-modules; or, from another point of view, the moduli stack $\mathcal{M}_{fmA}$ of formal $A$-modules. We ought to explain what this means. 
When $A$ is a commutative ring, a {\em formal $A$-module}
is a formal group law $F$ over a commutative $A$-algebra $R$, which is additionally equipped with a ring map $\rho: A \rightarrow \End(F)$ such that $\rho(a)(X) \equiv aX$ modulo $X^2$. An excellent introductory reference for formal $A$-modules is \cite{MR506881}. Higher-dimensional formal modules exist, but {\em all formal modules in this paper will be implicitly understood to be one-dimensional.}

Formal modules arise
in algebraic and arithmetic geometry, for example, in Lubin and Tate's famous theorem (in \cite{MR0172878}) on the abelian closure of a $p$-adic number field, in Drinfeld's generalizations of results of class field theory in \cite{MR0384707}, and in Drinfeld's $p$-adic symmetric domains, which are (rigid analytic) deformation
spaces of certain formal modules; see \cite{MR0422290} and \cite{MR1393439}.
Formal $A$-modules also arise in algebraic topology, by using the natural map from the moduli stack of formal $A$-modules to the moduli stack of formal groups to detect certain classes in the cohomology of the latter, particularly in order to resolve certain differentials in spectral sequences used to compute the Adams-Novikov $E_2$-term and stable homotopy groups of spheres; for example, see \cite{cmah7}.

More to the point for the present paper: it is easy to show (see \cite{MR0384707}) that there exists a classifying ring $L^A$ for formal $A$-modules, i.e., 
a commutative $A$-algebra $L^A$ such that $\hom_{A-alg}(L^A, R)$ is in natural bijection with the set of formal $A$-modules over $R$. It is not so easy to calculate $L^A$, however. Here is a summary of known results.
\begin{itemize}
\item The pioneer in this area was M. Lazard, who, in the case $A = \mathbb{Z}$, proved (see \cite{MR0393050}) that $L^{\mathbb{Z}} \cong \mathbb{Z}[x_1, x_2, \dots]$, a polynomial algebra on countably infinitely many generators. The ring $L^{\mathbb{Z}}$ is consequently often called the Lazard ring. 
\item Next, in \cite{MR0384707}, Drinfeld handled the case in which $A$ is the ring of integers in a local nonarchimedean field (e.g. a $p$-adic number field). In that case
Drinfeld proved that $L^A \cong A[x_1, x_2, \dots]$, again a polynomial algebra.
\item In \cite{MR506881}, Hazewinkel proved that the same result
holds for discrete valuation rings, as well as for global number rings of class number one. That is, for all such rings $A$, the classifying ring $L^A$ of formal $A$-modules is a polynomial $A$-algebra on countably infinitely many generators. 
\item Hazewinkel also makes the observation, in 21.3.3A of \cite{MR506881}, that the same result cannot possibly hold for
arbitrary global number rings. Specifically, when
$A$ is the ring of integers in 
$\mathbb{Q}(\sqrt[4]{-18})$, then Hazewinkel shows that
the sub-$A$-module of $L^A$ consisting of elements of grading degree $2$ (see Theorem \ref{basic existence thm on classifying hopf algebroids} for where this grading comes from) is not a free $A$-module, which could not occur if $L^A$ were polynomial. Hazewinkel does not, however,
attempt to compute $L^A$ for such rings $A$.
\item The preprint \cite{cmah2} contains calculations of $L^A$ for various Dedekind domains $A$, including cases with nontrivial class group. For a Dedekind domain $A$ of characteristic zero, it is shown that $L^A$ is always a symmetric $A$-algebra on a certain projective $A$-module, but fails to be a polynomial $A$-algebra when the relevant projective module is not free. 
\end{itemize}
The moduli stack $\mathcal{M}_{fmA}$ admits a natural presentation by $\Spec$ of the Hopf algebroid $(L^A, L^AB)$, and the flat cohomology of $\mathcal{M}_{fmA}$ coincides with the derived functors of the cotensor product (i.e., $\Cotor$) in the category of $L^AB$-comodules. Consequently some understanding of the ring $L^A$ and the Hopf algebroid $(L^A, L^AB)$ is a great help in calculating the cohomology of $\mathcal{M}_{fmA}$.

When formulated in terms of the Adams-Novikov spectral sequence, Adams's famous calculation of the image of the stable $J$-homomorphism, from \cite{MR0198470},  establishes that the flat cohomology group $H^1_{fl}(\mathcal{M}_{fm\mathbb{Z}}; \omega^n)$ is a finite abelian group of order equal to the denominator of the Riemann zeta-value $\zeta(1-n)$ times a power of $2$. Here $\omega$ is the line bundle of invariant differentials on $\mathcal{M}_{fm\mathbb{Z}}$. One would like to know if this result generalizes to rings other than $\mathbb{Z}$. 
In the paper \cite{MR745362}, Ravenel remarks that this result does ``not appear to generalize to other number fields. For example if the field is not totally real its Dedekind zeta function vanishes at all negative integers.'' 

This paper has two purposes:
\begin{enumerate}
\item to carry out a structural study of the Hopf algebroid $(L^A,L^AB)$ {\em for a general ring $A$}, aimed particularly at establishing structural properties of $(L^A,L^AB)$ which would support cohomology calculations, again {\em for a general ring $A$}, not just $\mathbb{Z}$, and not just a number ring or a Dedekind domain.
\item Subsequently, to use those structural properties, and some new computational tools, to make some cohomology calculations and demonstrate that the flat cohomology of the moduli of formal $A$-modules do in fact carry some zeta-function-theoretic data, {\em for a general ring $A$}.
\end{enumerate}

\subsubsection{Summary of structural results.}
\begin{description}
\item[Colimits] The functor sending $A$ to $L^A$ commutes with filtered colimits and with coequalizers, but not coproducts. The same is true for the functor sending $A$ to the Hopf algebroid $(L^A,L^AB)$. This is Proposition \ref{colimits for lazard rings}.
\item[Localization]  If $A$ is a commutative ring and $S$ a multiplicatively closed subset of $A$, then 
the homomorphism of graded rings $L^A[S^{-1}] \rightarrow L^{A[S^{-1}]}$ is an isomorphism. 
Furthermore, the homomorphism of graded Hopf algebroids
\begin{equation*}\label{map of hopf algebroids 22a} ( L^A[S^{-1}] , L^AB[S^{-1}]) \rightarrow (L^{A[S^{-1}]}, L^{A[S^{-1}]}B)\end{equation*}
is an isomorphism. This is Theorem \ref{lazard ring localization iso}. This particular result is not new: it appears also in Hazewinkel's book \cite{MR2987372}, although we think the proof we offer in this paper is a useful addition to the literature.
\item[Localization and cohomology] 
Let $A$ be a commutative ring and let $S$ be a multiplicatively closed subset of $A$. Let $f$ denote the stack homomorphism $f: \mathcal{M}_{fmA[S^{-1}]} \rightarrow \mathcal{M}_{fmA}$ classifying the underlying formal $A$-module of the universal formal $A[S^{-1}]$-module.
Then, for all quasicoherent $\mathcal{O}_{\mathcal{M}_{fmA}}$-modules $\mathcal{F}$, we have an isomorphism
\[ H^s_{fl}(\mathcal{M}_{fmA};  \omega^{\otimes t}\otimes \mathcal{F})[S^{-1}] \cong 
 H^s_{fl}(\mathcal{M}_{fmA[S^{-1}]};\omega^{\otimes t}\otimes  f^*\mathcal{F}) \]
for all integers $s,t$.
This is the stack-theoretic formulation of Corollary \ref{localization in cohomology}. It is also not new, having already appeared in \cite{MR745362} and in \cite{pearlmanthesis}.
Examples of its application in the course of cohomological calculations occur in Propositions \ref{2-loc prop} and \ref{torsion in LA2-buds}.
\item[Finite generation] 
\Cref{finiteness properties subsection} contains a variety of finiteness results which establish that, for a wide class of commutative rings $A$, the graded rings $L^A$ and $L^AB$ are finitely generated $A$-modules in each degree.
\item[Completion]  Theorem \ref{completion of the hopf algebroid} establishes that, for a wide class of commutative rings $A$, the functor sending $A$ to the Hopf algebroid $(L^A,L^AB)$ commutes with completion at a maximal ideal $I$.
As a consequence, Corollary \ref{ss completion corollary} establishes that the spectral sequence obtained from the $I$-adic filtration on the cobar complex of $(L^A,L^AB)$ converges to $H^*_{fl}(\mathcal{M}_{fm\hat{A}_I}; \mathcal{F})$, the cohomology of the moduli stack of formal $\hat{A}_I$-modules. The resulting spectral sequence is useful for making explicit calculations: see \cref{Construction of the 2-adic...} and \cref{cotor for Z} for examples.
\end{description}

\subsubsection{Summary of cohomological results.}

Extensive partial calculations of the flat cohomology of $\mathcal{M}_{fm\mathbb{Z}}$ have been made in stable homotopy theory, since \linebreak $H^*_{fl}(\mathcal{M}_{fm\mathbb{Z}}; \omega^{\otimes *})$ is the $E_2$-term of the Adams-Novikov spectral sequence which converges to the stable homotopy groups of spheres. Partial calculations of \linebreak $H^*_{fl}(\mathcal{M}_{fmA}; \omega^{\otimes *})$ for $A$ a (local or global) number ring can be found in \cite{MR745362}, \cite{formalmodules4}, and \cite{cmah7}, and for one particular number ring $A$, in \cite{MR2262886}. That is the extent of calculations of the cohomology of $\mathcal{M}_{fmA}$ to be found in the literature. 

In particular, there are no existing calculations of $H^*_{fl}(\mathcal{M}_{fmA}; \omega^{\otimes *})$ for rings $A$ of global dimension greater than one. There is good reason for this state of affairs: there are great difficulties\footnote{See the footnote at the start of \cref{Moduli...2-buds} for an indication of where these difficulties lie.} in calculating the ring $L^A$ when $A$ has global dimension $>1$. As a consequence, we begin our cohomological investigation of $\mathcal{M}_{fmA}$ with calculations of the flat cohomology of the moduli stack $\mathcal{M}_{fmA}^{2-buds}$ of {\em formal $A$-module $2$-buds}, i.e., power series in $R[[X,Y]]/(X,Y)^3$ which are required only to satisfy the formal group law axioms modulo $(X,Y)^3$, and which are equipped with an action of $A$ which again is only required to be unital and associative modulo $(X,Y)^3$. Compared to $\mathcal{M}_{fmA}$, the Artin stack $\mathcal{M}_{fmA}^{2-buds}$ is much simpler and easier to work with. Nevertheless it already enjoys remarkable cohomological properties in low degrees. Here are our results:
\begin{description}
\item[In degree zero] $H^0_{fl}(\mathcal{M}_{fmA}^{2-buds};\omega^{\otimes *})$, i.e., the sections of the tensor powers of the line bundle $\omega$ of invariant differentials, is described as follows. Let $A$ be a torsion-free commutative ring. In Theorem \ref{cotor0 thm} we obtain an isomorphism of $A$-modules:
\begin{align*}
 H^0_{fl}(\mathcal{M}_{fmA}^{2-buds}; \omega^{\otimes n})
  &\cong \left\{ \begin{array}{ll} 
   0 &\mbox{\ if\ } n<0 \\
   A &\mbox{\ if\ } n=0 \\
   0 &\mbox{\ if\ } n=1 \\
   \ker \tau_n^{I_2} 
    &\mbox{\ if\ } n>1 ,\end{array}\right. 
\end{align*} 
where $I_2$ is the ideal of $A$ generated by $2$ and by $a^2-a$ for all $a\in A$, and where $\tau_n^{I_2}$ is the $A$-module map
\begin{align*} 
 \tau_n^{I_2}: \Symm_A^{n}(I_2) &\rightarrow \Symm_A^{n-1}(I_2) \\
\nonumber x_1\cdot\dots\cdot x_n 
  &\mapsto \iota(x_1)x_2\cdot\dots\cdot x_n \\
\nonumber  &\ \ \ \ + x_1\cdot\iota(x_2)x_3\cdot\dots\cdot x_n \\
\nonumber  &\ \ \ \ + \dots + x_1\cdot x_2\cdot\dots \cdot \iota(x_{n-1})x_n.
\end{align*}
Here $\Symm_A^n(I_2)$ is the $n$th symmetric power of the ideal $I_2$, and $\iota$ denotes the inclusion of the ideal $I_2$ into $A$, regarded as an $A$-module morphism $\iota: I_2\rightarrow A$. See the discussion immediately preceding Theorem \ref{cotor0 thm} for some intuition behind the morphism $\tau_n^{I_2}$. See also \cref{Point-detecting ideals...} for discussion of the ideal $I_2$, including a universal property which it enjoys: it is the ``universal $\mathbb{F}_2$-point-detecting ideal,'' meaning it is the largest ideal of $A$ which is in the kernel of every ring homomorphism $A\rightarrow \mathbb{F}_2$.

Suppose $A$ is Noetherian. Then, in the particular case $n=2$, the kernel $\ker \tau_2^{I_2}$ coincides with the {\em delta-invariant} of the ideal $I_2$, which has enjoyed some attention in commutative algebra: see for example Micali-Roby in \cite{MR282964}, and Simis-Vasconcelos in \cite{MR610474}. The delta-invariant $\delta(I)$ is known (see \cite{MR992766}) to agree with the second Andre-Quillen homology group $H_2(A, A/I; A/I)$ of $A/I$ regarded as an $A$-algebra, with coefficients in $A/I$. 

Finitely generated ideals $I$ whose delta-invariant vanishes are called {\em syzygetic} in commutative algebra. We have Corollary \ref{H0 cor 2}: if $A$ is a Noetherian integral domain of characteristic zero, then $H^0_{fl}(\mathcal{M}_{fmA}^{2-buds}; \omega\otimes \omega)$ vanishes if and only if the universal $\mathbb{F}_2$-point-detecting ideal of $A$ is syzygetic.

If $A$ is a Cohen-Macaulay integral domain of characteristic zero, then Theorem \ref{H0 vanishing for CM} shows that $H^0_{fl}(\mathcal{M}_{fmA}^{2-buds}; \omega^{\otimes n})$ is trivial for all $n\neq 0$.
\item[The EFM spectral sequence] In Theorem \ref{change of endo-ring ss}, we construct a spectral sequence which converges to $H^*_{fl}(\mathcal{M}_{fmA}^{2-buds}; \omega^{\otimes *})$. Its $E_1$-term is a tensor product of $H^*_{fl}(\mathcal{M}_{fm\mathbb{Z}}^{2-buds}; \omega^{\otimes *})$ and symmetric powers of a certain module of twisted K\"{a}hler differentials in the sense of \cite{MR1975301}. Because this spectral sequence allows us to pass from the cohomology of the moduli of formal groups (i.e., formal $\mathbb{Z}$-modules) to the cohomology of the moduli of formal groups equipped with a larger ring (namely, $A$) of formal multiplications, we call it the ``extension of of formal multiplications spectral sequence,'' or ``EFM spectral sequence'' for short.
\item[Cohomology of certain twists] Assume that $A$ is a torsion-free commutative Noetherian ring. Using the EFM spectral sequence, in Theorem \ref{cotor1 calc thm} we calculate the flat cohomology in all degrees with coefficients in the first three tensor powers of $\omega$:
\begin{align*}
 H^*_{fl}(\mathcal{M}^{2-buds}_{fmA}; \omega^{\otimes n}) 
  &\cong 0\mbox{\ if\ } n<0,\\
 H^s_{fl}(\mathcal{M}^{2-buds}_{fmA}; \mathcal{O}) 
  &\cong \left\{ \begin{array}{ll} 
   A &\mbox{\ if\ } s=0 \\ 
   0 &\mbox{\ if\ } s\neq 0, \end{array}\right. \\
 H^s_{fl}(\mathcal{M}^{2-buds}_{fmA}; \omega) 
  &\cong \left\{ \begin{array}{ll} 
   A/I_2 &\mbox{\ if\ } s=1 \\ 
   0 &\mbox{\ if\ } s\neq 1, \end{array}\right. \\
 H^s_{fl}(\mathcal{M}^{2-buds}_{fmA}; \omega^{\otimes 2}) 
  &\cong \left\{ \begin{array}{ll} 
   \delta(I_2) &\mbox{\ if\ } s=0 \\ 
   A/I_2^2 &\mbox{\ if\ } s=1 \\ 
   A/I_2  &\mbox{\ if\ } s=2 \\
   0 &\mbox{\ otherwise},\end{array}\right. 
\end{align*}
Consequently we get Corollary \ref{local zeta-function cor}: $H^1_{fl}(\mathcal{M}_{fmA}^{2-buds}; \omega)$ is a finite abelian group of order equal to $2^{N_1}$, where $N_1$ is the number of $\mathbb{F}_2$-points of $\Spec A$, i.e., the logarithmic derivative of the $2$-local zeta-function $Z(\Spec A,t)$ of the affine scheme $\Spec A$, evaluated at $t=0$. See Remark \ref{remark on local zeta} for a very brief discussion of the local zeta-function of an affine variety.

The local zeta-function is, when evaluated at $t=p^{-s}$, an Euler factor in the Hasse-Weil zeta-function of a variety. The point here is that, even when restricting the scope of our calculations to the Artin stack of formal $A$-module $2$-buds, we still recover some zeta-function-theoretic information about $A$, not only for number rings $A$ (as Ravenel remarked about in \cite{MR745362}), but for characteristic zero integral domains quite generally. 

Some generalization to the number of $k$-points of $\Spec A$ for larger finite fields $k$ is possible, but requires calculations of the flat cohomology of $\mathcal{M}_{fmA}^{n-buds}$ for $n>2$, or the flat cohomology of $\mathcal{M}_{fmA}$. We have some preliminary results in this direction, but we regard them as beyond the scope of this already-too-long paper.
\item[Comparison to stable homotopy] The EFM spectral sequence calculates $H^*_{fl}(\mathcal{M}_{fmA}^{2-buds}; \omega^{\otimes *})$ in such a way that, along the way, it also calculates the homomorphism
\begin{align}\label{map in coh 0439}
 H^*_{fl}(\mathcal{M}_{fm\mathbb{Z}}^{2-buds}; \omega^{\otimes *}) &\rightarrow
 H^*_{fl}(\mathcal{M}_{fmA}^{2-buds}; \omega^{\otimes *})
\end{align}
induced by the map of Artin stacks $\mathcal{M}_{fmA}^{2-buds}\rightarrow\mathcal{M}_{fm\mathbb{Z}}^{2-buds}$ classifying the underlying formal group law $2$-bud of the universal formal $A$-module $2$-bud. There is also a stack map $\mathcal{M}_{fm\mathbb{Z}}\rightarrow \mathcal{M}_{fm\mathbb{Z}}^{2-buds}$ classifying the underlying $2$-bud of the universal formal group law. The flat cohomology of $\mathcal{M}_{fm\mathbb{Z}}$ is the input for the Adams-Novikov spectral sequence converging to the stable homotopy groups of spheres. In \cref{COER ss} we calculate that the element $\eta$ in $H^1_{fl}(\mathcal{M}_{fm\mathbb{Z}}; \omega)$ which yields the generator of the first stable homotopy group $\mathbb{Z}/2\mathbb{Z}\in \pi_1^{st}(S^0)$ is the image, under the map 
\begin{align*}
 H^1_{fl}(\mathcal{M}_{fm\mathbb{Z}}^{2-buds}; \omega) &\rightarrow H^1_{fl}(\mathcal{M}_{fm\mathbb{Z}}; \omega), \end{align*}
of a unique element $H^1_{fl}(\mathcal{M}_{fm\mathbb{Z}}^{2-buds}; \omega)$.
The image of this element in \linebreak $H^1_{fl}(\mathcal{M}_{fmA}^{2-buds}; \omega)$
is an $A$-module generator. The $A$-module $H^1_{fl}(\mathcal{M}_{fm\mathbb{A}}^{2-buds}; \omega)$ has order equal to $2^{N_1}$, as already explained. In this sense, the same element in the flat cohomology of $\mathcal{M}_{fm\mathbb{Z}}$ which is responsible for the first stable homotopy group of spheres is also responsible for the first coefficient in the $2$-local zeta function for any torsion-free commutative Noetherian ring.
\end{description}

The author of this paper is an algebraic topologist, but is optimistic that this paper may be interesting to readers from various mathematical backgrounds. I have made an effort to write this paper so that it will be readable by researchers in other subjects. However, in computational algebraic topology, there is a tradition of engaging in detailed and meticulous spectral sequence calculations, and to some extent that tradition is reflected in the second half of this paper. I apologize to readers who do not care for spectral sequence calculations, and I am grateful for their patience with \cref{Cohomology calculations.} and \cref{H0 and H1 of the moduli of...}, and with the length of this paper.

\begin{conventions}\leavevmode
\begin{itemize}
\item We often deal with graded rings in this paper, but we will always set up the gradings on the classifying rings $L^A$ and on $L^AB$ so that all elements are in even degrees. Consequently the Koszul graded-commutativity sign relation does not occur
in our discussions of the structure of $L^A$ and $L^AB$. We often use the phrase ``commutative graded ring''  rather than ``graded-commutative'' to emphasize that the graded ring in question is assumed to be {\em strictly} commutative, not merely commutative up to sign.
\item We call a ring {\em torsion-free} when its underlying abelian group is torsion-free.
\end{itemize}
\end{conventions}

\subsection{Review of standard facts about $L^A$ and $L^AB$.}\label{review subsection}

Nothing in this subsection is new, but we think it may be helpful to the reader to have many of the basic ideas and known results on formal modules and their classifying Hopf algebroids collected in one place.

\subsubsection{Formal modules and formal module $n$-buds.}
\label{Formal modules and formal module n-buds.}
If $A$ is a commutative ring and $R$ is a commutative $A$-algebra,
then a (one-dimensional) {\em formal $A$-module\footnote{Morally, a formal $A$-module $F$ is a ``formal group law with complex multiplication by $A$.'' This perspective was taken already by Lubin and Tate in \cite{MR0172878}, and suggests the close connection between formal modules and abelian varieties with particular endomorphism rings. The turn of phrase ``formal group law with complex multiplication'' is not as clear as might be hoped, however: if $A$ is the ring of integers in a totally real extension of $\mathbb{Q}$, then an abelian variety with a suitable action by $A$ would be called an abelian variety with {\em real} multiplication, rather than {\em complex} multiplication. Consequently we prefer to write ``formal module'' rather than ``formal group law with complex muliplication.'' Rather than complex multiplication or real multiplication, we will refer to the action of $A$ on $F$ as {\em formal multiplication}.} over $R$} is a formal group law $F$ over $R$, together with 
a ring homomorphism $\rho: A \rightarrow \End(F)$
such that $\rho(a)\in \End(F) \subseteq R[[X]]$
is congruent to $aX$ modulo $X^2$. 
Here $\End(F)$ is a ring in which:
\begin{itemize}
\item addition is given by formal addition, i.e., the sum of $f(X)$ and $g(X)$ is $F(f(X), g(X))$, {\em not} the 
ordinary componentwise addition of power series, 
\item and multiplication is given by composition of power series, {\em not} the usual multiplication of power series in 
$R[[X]]$.
\end{itemize}

If $n$ is a positive integer, a {\em formal $A$-module $n$-bud over $R$}
is a formal group law $n$-bud over $R$, i.e., an element $F(X,Y)\in R[[X,Y]]/(X,Y)^{n+1}$ which satisfies the unitality, associativity, commutativity, and existence of inverses axioms modulo $(X,Y)^{n+1}$,
together with a ring homomorphism $\rho: A \rightarrow \End(F)$
such that the endomorphism $\rho(a)\in \End(F) \subseteq R[[X]]/(X^{n+1})$
is congruent to $aX$ modulo $X^2$. 

In this paper we will {\em always} write $\End(F)$ for the endomorphism ring of a formal group law $F$ or formal group law $n$-bud $F$. That is, even if $F$ has the additional structure of a formal module, by $\End(F)$ we will mean the endomorphism ring of
$F$ as a formal group law or formal group law $n$-bud, without regard to any additional structure.

\subsubsection{Hopf algebroids and stacks.}
\label{Hopf algebroids and stacks.}
This paper is largely about certain graded Hopf algebroids, i.e.,
cogroupoid objects in commutative graded rings. 
We give only a cursory review of the basic theory here. For more detail, we refer readers to the standard reference for Hopf algebroids and their cohomology, Appendix 1 of \cite{MR860042}. Whenever convenient, we will use the common notations for structure maps of bialgebroids and Hopf algebroids: $\eta_L$ for left unit, $\eta_R$ for right unit, $\Delta$ for coproduct, and $\epsilon$ for augmentation. 

Starting in \cref{Calculation of Cotor zero}, we will make calculations of cohomology groups of Hopf algebroids. The cohomology groups of a graded Hopf algebroid $(A,\Gamma)$ with coefficients in a graded $\Gamma$-comodule $M$ can be defined in two (isomorphic) ways:
\begin{itemize}
\item as the right derived functors $\Cotor^{s,t}_{\Gamma}(A,M) = \Cotor^{s}_{\Gamma}(A,\Sigma^t M)$ of the cotensor product $A\Box_{\Gamma}-: \gr\Comod(\Gamma)\rightarrow \Ab$ applied to $M$,
\item or as the relative right derived functors $\Ext^{s,t}_{(A,\Gamma)}(A,M) = \Ext^{s}_{(A,\Gamma)}(A,\Sigma^t M)$ of $\hom_{\gr\Comod(\Gamma)}(A,-): \gr\Comod(\Gamma)\rightarrow \Ab$ applied to $M$,
relative to the allowable class generated by the comodules tensored up from $A$. This is a relative 
$\Ext$-group, in the sense of relative homological algebra, as in Chapter IX of \cite{MR1344215},
\item A third description of the cohomology of $(A,\Gamma)$ with coefficients in $M$ is available whenever the unit maps $\eta_L,\eta_R$ of the Hopf algebroid are smooth (respectively, formally smooth). In that case, the stackification $\mathcal{X}$ of the groupoid scheme $(\Spec A,\Spec \Gamma)$ is an Artin (respectively, formally Artin) stack in the fpqc topology. Its category of quasicoherent modules is equivalent to the category of $\Gamma$-comodules. Under this equivalence, the Hopf algebroid cohomology group $\Cotor^s_{\Gamma}(A,M)$ is isomorphic to the flat stack cohomology group $H^s_{fl}(\mathcal{X}; \tilde{M})$, where $\tilde{M}$ is the quasicoherent module associated to the comodule $M$. These facts are standard; a nice reference is \cite{pribblethesis}.
\end{itemize}

We also refer to ``the moduli stack of formal $A$-modules'' several times
in this paper. This is slightly ambiguous for the following reason:
formal $A$-modules have only a moduli prestack and not a moduli stack.
This moduli prestack ``stackifies'' (as in \cite{MR1771927})
to a stack which is a moduli stack for ``coordinate-free'' formal $A$-modules\footnote{Here is a bit of detail about what a ``coordinate-free'' formal module is. These details are routine, not important for the rest of this paper, and can be safely skipped by the reader.

It is classical that a formal group {\em law} over $R$ is a power series $F(X,Y)\in R[[X,Y]]$ satisfying associativity, commutativity, unitality, and inverse axioms, i.e., $F(X,Y)$ defines the structure of a commutative group object on the affine formal scheme $\Spf R[[X]]$, with identity element $0$. A {\em formal group} is a commutative group structure on the formal affine line $\hat{\mathbb{A}}_R^1$ with identity element $0$. Consequently a formal group law is a formal group together with a choice of isomorphism $\hat{\mathbb{A}}_R^1\cong\Spf R[[X]]$. Formal groups form a stack $\mathcal{M}_{fg}$, while formal group {\em laws} only form a prestack whose stackification is equivalent to $\mathcal{M}_{fg}$.

A similar story applies here. It is usual to say ``formal $A$-module'' to mean the power series $F(X,Y)\in R[[X,Y]]$ equipped with the action of $A$ by further power series over $R$, as defined in \cref{Formal modules and formal module n-buds.}. Perhaps it would better if such objects were instead called ``formal $A$-module laws,'' since such an object determines the structure of an $A$-module object on the affine formal scheme $\Spf R[[X]]$, with identity element $0$. Then we could use the term ``formal $A$-module'' to mean the structure of an $A$-module object on $\hat{\mathbb{A}}^1_R$ with identity element $0$, and the terminology would mirror the standard distinction between ``formal group law'' and ``formal group.'' Using the terms in this way, formal $A$-module laws form a prestack, represented by the groupoid scheme $(\Spec L^A, \Spec L^AB)$. The stackification of that prestack is then equivalent to the stack of formal $A$-modules $\mathcal{M}_{fmA}$. Unfortunately the weight of tradition is against this distinction in terminology between ``formal modules'' and ``formal module laws.''},
a situation which parallels that of formal group laws and formal groups, as in \cite{smithlingthesis}.

\subsubsection{The Hopf algebroid $(L^A,L^AB)$.}

Theorem \ref{basic existence thm on classifying hopf algebroids}
is the main foundational result about the Hopf algebroid $(L^A,L^AB)$. It
gathers together many results proven in chapter 21 of \cite{MR506881},
although parts of the theorem are older than Hazewinkel's
book; for example, the computation of the ring
$L^A$, when $A$ is a field or the ring of integers in a nonarchimedean local field, is due to Drinfeld in \cite{MR0384707}.
\begin{theorem}\label{basic existence thm on classifying hopf algebroids}
Let $A$ be a commutative ring.
\begin{itemize}
\item
Then there exist commutative $A$-algebras $L^A$ and $L^AB$
having the following properties:
\begin{itemize}
\item For any commutative $A$-algebra $R$, there exists
a bijection, natural in $R$, between the set 
of $A$-algebra homomorphisms $L^A\rightarrow R$ and
the set of formal $A$-modules over $R$.
\item For any commutative $A$-algebra $R$, there exists
a bijection, natural in $R$, between the set 
of $A$-algebra homomorphisms $L^AB\rightarrow R$ and
the set of strict\footnote{Recall that an isomorphism $f(X)$ of formal groups, or of formal $A$-modules, is said to be {\em strict} if $f(X) \equiv X\mod X^2.$} isomorphisms of formal $A$-modules over $R$.
\end{itemize}
\item The natural maps of sets between the set of
formal $A$-modules over $R$ and the set of strict isomorphisms of 
formal $A$-modules over $R$ (sending a strict isomorphism to 
its domain or codomain; or sending a formal module to its identity
strict isomorphism; or composing two strict isomorphisms; or
sending a strict isomorphism to its inverse) are
co-represented by maps of $A$-algebras between $L^A$ and $L^AB$.
Consequently $(L^A, L^AB)$ is a Hopf algebroid
co-representing the functor sending a commutative $A$-algebra $R$
to its groupoid of formal $A$-modules and their strict isomorphisms.
\item If $n$ is a positive integer, then the functor from commutative $A$-algebras to groupoids which sends a commutative $A$-algebra $R$ to 
the groupoid of formal $A$-module $n$-buds over $R$ and strict isomorphisms is also co-representable by a Hopf algebroid
$(L^A_{n-buds}, L^A_{n-buds}B)$. Since the groupoid of formal $A$-modules over $R$ is the inverse limit over $n$ of the 
groupoid of formal $A$-module $n$-buds over $R$, we have that
\[ (L^A, L^AB) \cong \left(\colim_{n\rightarrow \infty} L^A_{n-buds},\colim_{n\rightarrow \infty} L^A_{n-buds}B\right).\]
For example, $L^A_{\leq 1} \cong A$ as commutative $A$-algebras.

The filtration of $L^A$ and $L^AB$ by $L^A_{n-buds}$ and $L^A_{n-buds}B$ induces a grading on $L^A$ and on $L^AB$, in which
the {\em indecomposable} homogeneous grading degree $2n$ elements in $L^A$ are the parameters for deforming (i.e., extending) a formal $A$-module $n$-bud to
a formal $A$-module $(n+1)$-bud. 
The summands of $L^A$ and of $L^AB$ of odd grading degree are trivial.
\item 
If $A$ is a field of characteristic zero or a discrete valuation ring or a global number ring of class number one,
then we have isomorphisms of graded $A$-algebras
\begin{align*} 
 L^{A}_{n-buds}
  &\cong {A}[x_1^{A},x_2^{A},x_3^{A},\dots , x_n^A], \\
 L^{A}_{n-buds}B
  &\cong L^{A}_{n-buds}[t_1^{A},t_2^{A},t_3^{A},\dots ,t_n^A] , \mbox{\ \ and\ consequently} \\
 L^{A} 
  &\cong {A}[x_1^{A},x_2^{A},x_3^{A},\dots ], \\
 L^{A}B
  &\cong L^{A}[t_1^{A},t_2^{A},t_3^{A},\dots ] ,
\end{align*}
with each $x_i^A$ and each $t_i^A$ homogeneous of grading degree $2i$.
(However, the natural map $L^A\rightarrow L^B$ induced by a ring homomorphism $A\rightarrow B$ does not necessarily send each $x_i^A$ to $x_i^B$!)
\end{itemize}
\end{theorem}
The factor of 2 in the gradings in Theorem \ref{basic existence thm on classifying hopf algebroids} is due to the graded-commutativity
sign convention in algebraic topology and the fact that
$L^{\mathbb{Z}}$, with the above grading, 
is isomorphic to the graded ring of 
homotopy groups $\pi_*(MU)$ of the complex bordism spectrum $MU$,
while $L^{\mathbb{Z}}B$ with the above grading is isomorphic to the graded ring
$\pi_*(MU\smash MU)$ of stable co-operations in complex bordism. In fact
$(L^{\mathbb{Z}}, L^{\mathbb{Z}}B) \cong (\pi_*(MU), \pi_*(MU \smash MU))$ as graded Hopf algebroids.
See \cite{MR0253350} for these ideas. In the base case $A=\mathbb{Z}$, one often writes $L$ and $LB$ rather than $L^{\mathbb{Z}}$ and $L^{\mathbb{Z}}B$.

Proposition \ref{drinfeld presentation} appears as Proposition 1.1 in \cite{MR0384707}.
\begin{prop}\label{drinfeld presentation} Let $A$ be a commutative ring, let $n$ be an integer, and let $D^A$ denote the homogeneous ideal in $L^A$ generated by all products of elements $xy$ with $x,y\in L^A$ each homogeneous of positive degree.
Let $\overline{L}^A$ denote the quotient ring $L^A/D^A$. 
The ring $\overline{L}^A$ is graded, so we may consider its degree $n$ summand $\overline{L}^A_n$ for various integers $n$.
The ring $L^A$ is concentrated in even degrees, so $\overline{L}^A$ is as well.
If $n\geq 2$, then $\overline{L}^A_{2n-2}$ is isomorphic
to the $A$-module generated by symbols $\gamma$
and $\{ c_a: a\in A\}$, that is, one generator $c_a$ for each element
$a$ of $A$ along with one additional generator $\gamma$,
modulo the relations:
\begin{align}
\label{hazewinkel relation 20} (a^n - a)\gamma &= \nu(n)c_a \mbox{\ \ for\ all\ } a\in A \\
\label{hazewinkel relation 21} c_{a+b}-c_a-c_b &= \gamma\frac{(a+b)^n - a^n - b^n}{\nu(n)} \mbox{\ \ for\ all\ } a,b\in A \\
\label{hazewinkel relation 22} ac_b + b^nc_a &= c_{ab} \mbox{\ \ for\ all\ } a,b\in A ,\end{align}
where $\nu(n)$ is defined to be the integer $1$ if $n$ is not a prime power, while $\nu(n)=p$ if $n$ is a power of a prime number $p$.
We will call this
{\em Drinfeld's presentation for $\overline{L}^A_{2n-2}$.}\end{prop}

The grading degrees in Proposition \ref{drinfeld presentation} are twice what they are in Drinfeld's statement
of the result in \cite{MR0384707}. Our gradings are chosen to match the gradings that occur in algebraic topology, where the generator of $L^{\mathbb{Z}} \cong MU_*$ classifying an extension of a formal group $n$-bud to a formal group $(n+1)$-bud is in grading degree $2n$ rather than $n$.

\subsection{Change of $A$.}

Proposition \ref{tensoring on one side} is a standard tool in Hopf algebroids (see A1.3.12 of \cite{MR860042}), and we omit the proof. 
\begin{prop} \label{tensoring on one side} Let $(R,\Gamma)$ be a commutative bialgebroid over a commutative ring $A$, and let $S$ be a right
$\Gamma$-comodule algebra, such that the following diagram commutes:
\begin{equation}  \label{right comodule alg extends unit} \xymatrix{ R\ar[r]^{\eta_R}\ar[d]^f & \Gamma\ar[d]^{f\otimes_R \id_\Gamma} \\ S\ar[r]^\psi & S\otimes_R \Gamma}\end{equation}
where $f$ is the $R$-algebra structure map $R\stackrel{f}{\longrightarrow} S$. 
Then we have a bialgebroid $(S,S\otimes_R\Gamma)$, with right unit $S\rightarrow S\otimes_R\Gamma$ 
equal to the comodule structure map $\psi$ on $S$. The map 
\begin{equation}\label{bialgebroid map 1} (R, \Gamma) \rightarrow (S, S\otimes_R \Gamma),\end{equation}
with components $f$ and $\psi$, is a morphism of bialgebroids.

 If $(R,\Gamma)$ is a Hopf algebroid (respectively, graded Hopf algebroid), then so is $(S,S\otimes_R\Gamma)$, and \eqref{bialgebroid map 1} is a map of Hopf algebroids (respectively, graded Hopf algebroids).
%

If, furthermore, the following conditions are also satisfied:
\begin{itemize}
\item $(R,\Gamma)$ is a graded Hopf algebroid which is connected
(i.e., the grading degree zero summand $\Gamma^0$ of $\Gamma$
is exactly the image of $\eta_L : R \rightarrow \Gamma$,
equivalently $\eta_R: R \rightarrow \Gamma$), and
\item $S$ is a graded $R$-module concentrated in degree zero, and
\item $N$ is a graded left $S\otimes_R\Gamma$-comodule which is 
flat as an $S$-module, and
\item $M$ is a graded right $\Gamma$-comodule, 
\end{itemize}
then we have an isomorphism
\[ \Ext^{s,t}_{(R, \Gamma)}(M, N) \cong \Ext^{s,t}_{(S, S\otimes_R \Gamma)}(S\otimes_R M, N)\]
for all nonnegative integers $s$ and all integers $t$.
\end{prop}

Proposition \ref{formal module base change thm} appeared originally in the unpublished doctoral thesis \cite{pearlmanthesis} of A. Pearlman:
\begin{prop}\label{formal module base change thm}
Let $f: A \rightarrow A^{\prime}$ be a homomorphism of commutative rings.
Then $L^{A^{\prime}}$ admits a canonical $L^AB$-comodule structure satisfying the conditions of Proposition \ref{tensoring on one side}. Consequently we have an isomorphism of graded Hopf algebroids 
\[ (L^{A^{\prime}}, L^{A^{\prime}}B) \cong (L^{A^{\prime}}, L^{A^{\prime}} \otimes_{L^A} L^AB),\]
and an isomorphism in cohomology
\[ \Cotor^{s,t}_{L^AB}(M, N) \cong \Cotor^{s,t}_{L^{A^{\prime}}B}(M\otimes_{L^A} L^{A^{\prime}}, N)\]
for all nonnegative integers $s$, all integers $t$, any graded right $L^AB$-comodule $M$, and any graded $L^{A^{\prime}}B$-comodule $N$ which is flat as 
a $L^{A^{\prime}}$-module.
\end{prop}

\section{Generalities on $L^A$ and $L^AB$.}\label{generalities on LA and LAB}

\subsection{Colimits.}

\begin{prop}\label{colimits for lazard rings}
Let $\mathcal{L},\mathcal{LB}$ be the functors
\begin{align*}
 \mathcal{L}: \Comm\Rings &\rightarrow \Comm\Rings \\
  \mathcal{L}(A) &= L^A \\
 \mathcal{LB}: \Comm\Rings &\rightarrow \Comm\Rings \\
  \mathcal{LB}(A) &= L^AB .
\end{align*}
Then $\mathcal{L}$ and $\mathcal{LB}$ each commute with filtered colimits, and $\mathcal{L}$ and $\mathcal{LB}$ each commute with 
coequalizers.
\end{prop}
\begin{proof}
Let $\mathcal{D}$ be a small category. Suppose that either $\mathcal{D}$ is filtered or $\mathcal{D}$ is the category indexing a parallel pair,
i.e., the Kronecker quiver 
\[ \xymatrix{ \bullet \ar@<1ex>[r] \ar@<-1ex>[r] & \bullet .}\]
Let $G: \mathcal{D}\rightarrow\Comm\Rings$ be a functor, let $R$ be a commutative ring, 
and suppose we are
given a cone $\mathcal{L} \circ G \rightarrow R$.
Then $R$ has the natural structure of a commutative $\colim G$-algebra, since the grading degree zero 
subring of each $\mathcal{L}(G(d))$ is isomorphic to the ring $G(d)$ itself.
Since $\mathbb{Z}$ is initial in commutative rings, there is a unique cocone $\mathbb{Z}\rightarrow G$
and hence a canonical cocone $L^{\mathbb{Z}} \rightarrow \mathcal{L}\circ G$.
Hence the cone $\mathcal{L} \circ G \rightarrow R$ describes a choice of formal
group law $F$ over the commutative $\colim G$-algebra $R$, together with a choice of 
ring map $\rho_d: G(d)\rightarrow \End(F)$ for each $d\in \ob\mathcal{D}$, compatible with the
morphisms in $\mathcal{D}$, and such that $\rho_d(r)(X) \equiv rX$ modulo $(X^2) \subseteq (\colim G)[[X]]$ for all $r\in G(d)$.

The colimit $\colim G$ here is computed in commutative rings. However, the ring $\End(F)$ is typically not commutative, so the universal property of $\colim G$ does not automatically yield a ring map $\colim G \rightarrow \End(F)$. We need one extra step before we get such a ring map: we observe that the image $\im\rho_d$ of each $\rho_d$ is a commutative subring of $\End(F)$, so the union of the family of subrings
$\cup_{d\in \ob\mathcal{D}} \im \rho_d$ is a commutative subring of $\End(F)$ since $\mathcal{D}$ is either filtered or is the category indexing
parallel pairs\footnote{This is the only part of the argument that uses the assumption that $\mathcal{D}$ is either filtered or the Kronecker quiver. The argument fails if $\mathcal{D}$ is only assumed to be an arbitrary small category. As far as I know there is no reason to believe that the conclusion of 
Proposition \ref{colimits for lazard rings} holds for coproducts, precisely because of the distinction between coproducts
in commutative rings and coproducts in associative rings.}.
Hence we have a cone $G \rightarrow \cup_{d\in \ob\mathcal{D}}\im\rho_d$ in the category of commutative rings, hence a canonical map
$\rho: \colim G \rightarrow \cup_{d\in \ob\mathcal{D}}\im\rho_d$ such that $\rho(r)(X) \equiv rX$ modulo $X^2$ for all $r\in \colim G$,
hence $F$ is a formal $\colim G$-module over $R$.
Clearly if we began instead with a formal $\colim G$-module over $R$, by neglect of structure we get
a cone $\mathcal{L}\circ G\rightarrow R$, and the two operations (sending such a cone to its $\colim G$-module, and sending the $\colim G$-module to its cone) are mutually inverse. So $\colim(\mathcal{L}\circ G) \cong \mathcal{L}(\colim G)$.

For $\mathcal{LB}$: we have already seen in Proposition \ref{formal module base change thm} that $\mathcal{LB}$ is naturally equivalent
to the functor $\mathcal{L} \otimes_{L^{\mathbb{Z}}} L^{\mathbb{Z}}B$. Since base change commutes with arbitrary colimits of commutative rings,
the fact that $\mathcal{L}$ commutes with filtered colimits and coequalizers implies the same for $\mathcal{LB}$.
\end{proof}

\begin{remark}
Proposition \ref{colimits for lazard rings} provides, at least in principle, a means of computing $L^A$ and $L^AB$ for all
commutative rings $A$: first, represent $A$ as the coequalizer of a pair of maps
\begin{equation}\label{}\xymatrix{
\mathbb{Z}[G]
   & 
  \mathbb{Z}[R] \ar@<1ex>[l]\ar@<-1ex>[l]
    }\end{equation}
where $G$ is a set of generators and $R$ a set of relations, and $\mathbb{Z}[G],\mathbb{Z}[R]$ are the free commutative algebras 
generated by the sets $G$ and $R$, respectively. Then $L^A$ is just the coequalizer, in commutative rings, of the two
resulting maps $L^{\mathbb{Z}[R]}\rightarrow L^{\mathbb{Z}[G]}$. 

Consequently, if one can compute $L^{A}$ for polynomial rings $A$, then one can (at least in principle) compute $L^A$ for all
commutative rings $A$. Unfortunately, the computation of $L^A$ for polynomial rings $A$ is quite difficult, 
and since the functor $\mathcal{L}$ does not commute with coproducts, it is not as simple as computing
$L^{\mathbb{Z}[x]}$ and then taking an $n$-fold tensor power to get $L^{\mathbb{Z}[x_1, \dots ,x_n]}$.
\end{remark}

\subsection{Localization.}

In the proof of Theorem 21.3.5 of Hazewinkel's excellent book \cite{MR506881}, also appearing in the second edition \cite{MR2987372}, 
one finds the following statement:
\begin{quotation}
  ``By the very definition of $L_A$ (as the solution of a certain universal problem) we have that $(L_A)_{\mathfrak{p}} = L_{A_{\mathfrak{p}}}$
for all prime ideals $\mathfrak{p}$ of $A$.''
\end{quotation}
It is true that $(L_A)_{\mathfrak{p}}$ is isomorphic to $L_{A_{\mathfrak{p}}}$. However, we prefer to give a few more lines of proof, since the universal properties of these rings do not obviously imply that every formal $A$-module over a commutative $A_{\mathfrak{p}}$-algebra extends to a formal $A_{\mathfrak{p}}$-module, as the
endomorphism ring $\End(F)$ of a formal group law defined over a ring $R$ is typically not an $R$-algebra. This is clear from the famous example of the endomorphism ring of a height $n$ formal group law over $\mathbb{F}_{p^n}$ being the maximal order in the invariant $1/n$ central division algebra over $\mathbb{Q}_p$, which is certainly not an $\mathbb{F}_{p^n}$-algebra. 
The claimed isomorphism also does not follow from Proposition \ref{colimits for lazard rings}, the fact that $A\mapsto L^A$ commutes with coequalizers and filtered colimits, since although localizations of {\em modules} can be defined as colimits in that category of modules, a localization of a {\em commutative ring} is not usually expressible as a colimit in the category of commutative rings. The morphisms in the diagram whose colimit computes the localization of the underlying module typically fail to be ring homomorphisms.

However, suppose that $r$ is a unit in a commutative $A$-algebra $R$, and suppose that $F$ is a formal $A$-module over $R$. Then the power series $\rho(r)(X)\in R[[X]]$ admits a (unique) composition inverse, by the tangent condition $\rho(r)(X) \equiv rX \mod X^2$ on the formal group law. Consequently, if $S$ is a set of non-zero-divisors in $A$, and if we write $G$ for the universal formal $A$-module base-changed to $L^A\otimes_A A[S^{-1}]$, then the formal $A$-multiplication map $\rho: A \rightarrow \End(G)$ admits a unique extension to a ring homomorphism $A[S^{-1}]\rightarrow \End(G)$, which furthermore satisfies the tangent condition. Hence $G$ is in fact a formal $A[S^{-1}]$-module. Consequently every formal $A$-module over a commutative $A[S^{-1}]$-algebra admits a unique compatible formal $A[S^{-1}]$-multiplication. From here it is routine to see how Hazewinkel's argument establishes the isomorphism $(L^A)[S^{-1}] \cong L^{A[S^{-1}]}$, and consequently the theorem:
\begin{theorem}\label{lazard ring localization iso}
Let $A$ be a commutative ring and let $S$ be a multiplicatively closed subset of $A$. 
Then the homomorphism of graded rings $L^A[S^{-1}] \rightarrow L^{A[S^{-1}]}$ is an isomorphism. 
Even better, the homomorphism of graded Hopf algebroids
\begin{equation*}
( L^A[S^{-1}] , L^AB[S^{-1}]) \rightarrow (L^{A[S^{-1}]}, L^{A[S^{-1}]}B)\end{equation*}
is an isomorphism of Hopf algebroids.
\end{theorem}

\begin{corollary}\label{localization in cohomology}
Let $A$ be a commutative ring and let $S$ be a multiplicatively closed subset of $A$. 
Then, for all graded left $L^A[S^{-1}]$-comodules $M$, we have an isomorphism
\[ \left(\Cotor^{s,t}_{L^AB}(L^A, M)\right)[S^{-1}] \cong \Cotor^{s,t}_{L^{A[S^{-1}]}B}(L^{A[S^{-1}]}, M)\]
for all nonnegative integers $s$ and all integers $t$.
\end{corollary}

\subsection{Finiteness, separation, and completion properties.}\label{finiteness properties subsection}

\begin{lemma}\label{finite generation lemma}
Let $A$ be a commutative ring, and let $R$ be a commutative graded $A$-algebra which is connective, i.e., the degree $n$ grading summand
$R_n$ is trivial for all $n<0$. Suppose that, for all integers $n$, the $A$-module $R_n/D_n$ is finitely generated, where
$D_n$ is the sub-$A$-module of $R_n$ generated by all elements of the form $xy$ where $x,y$ are homogeneous elements of $R$ of 
grading degree $<n$.

Then, for all integers $n$, $R_n$ is a finitely generated $A$-module.
\end{lemma}
\begin{proof}
Routine.
\end{proof}

\begin{prop}\label{finite typeness of lazard ring}
Let $A$ be a torsion-free commutative ring, and suppose that
$A$ is finitely generated as a commutative ring. 
Then, for each integer $m$, the degree $m$ summand $L^A_m$
of the classifying ring $L^A$ of formal $A$-modules is
a finitely generated $A$-module.
\end{prop}
\begin{proof}
Suppose that $A$ is generated, as a commutative ring, by a finite set of 
generators $x_1, \dots ,x_n$. Using Drinfeld's presentation for $\overline{L}^A_{2m-2}$, 
the Drinfeld relations \eqref{hazewinkel relation 21} and \eqref{hazewinkel relation 22} yield that $\overline{L}^A_{2m-2}$ is generated, as an $A$-module,
by the $n+1$ elements $\gamma, c_{x_1}, \dots ,c_{x_n}$.
(The assumption that $A$ is torsion-free is being used with relation \eqref{hazewinkel relation 21}, so that
there is at most one way of dividing $((a+b)^m - a^m - b^m)\gamma$ by $\nu(n)\in\mathbb{Z}$, so that $c_{a+b}$ is {\em uniquely} determined
by $\gamma$ and $c_{a}$ and $c_{b}$.)
Now Lemma \ref{finite generation lemma} implies that $L^A_m$ is a finitely generated $A$-module for all integers $m$.
\end{proof}
In Proposition \ref{finite typeness of lazard ring} it is important
that $L^A$ is typically not a finitely-generated $A$-module,
nor even finitely generated as an $A$-algebra; rather,
the summand in each individual degree is a finitely generated
$A$-module. 

\begin{corollary}\label{finite generation and separatedness of local lazard ring}
Let $A$ be a torsion-free commutative ring, and suppose that
$A$ is finitely generated as a commutative ring.
Let $I$ be a maximal ideal of $A$, and let $A_I$ denote
$A$ localized at $I$, i.e., $A$ with all elements outside of 
$I$ inverted.
Then, for each integer $m$, the grading degree $m$ summand
$L^{A_I}_m$ of the classifying ring $L^{A_I}$ of formal $A_I$-modules is
a finitely-generated, 
$I$-adically separated $A_I$-module.
\end{corollary}
\begin{proof}
By Propositions \ref{lazard ring localization iso} and \ref{finite typeness of lazard ring}, $(L^A_m)_I \cong L^{A_I}_m$ is a finitely generated $A_I$-module
for all integers $m$.
The ring $A_I$ is Noetherian and local, so the Krull intersection theorem (classical; see Corollary 10.20 in \cite{MR0242802}) implies that
every finitely generated $A_I$-module is $I$-adically separated.
\end{proof}

\begin{prop}\label{hensel local ring lazard ring finite generation}
Let $A$ be a torsion-free Henselian local commutative ring with maximal ideal $\mathfrak{m}$.
Suppose that $\mathfrak{m}$ can be generated by $\kappa$ elements, where $\kappa$ is some cardinal number.

Then, for each positive integer $n$, the $A$-module 
$\overline{L}^A_{2n-2}$ can be generated by:
\begin{itemize}
\item $1+\kappa$ elements, if the residue field $A/\mathfrak{m}$ is isomorphic to a finite field $\mathbb{F}_q$ and $n$ is a power of $q$,
\item and $1$ element (i.e., $\overline{L}^A_{2n-2}$ is a cyclic $A$-module) otherwise.
\end{itemize}
\end{prop}
\begin{proof}
For this theorem we use Drinfeld's presentation for $\overline{L}^A_{2n-2}$. 
Let $p$ denote the characteristic of $A/\mathfrak{m}$. (We allow $p=0$ as a possibility.)
There are three cases to consider:
\begin{itemize}
\item {\bf If $n$ is not a power of $p$:} Then $\nu(n)$ is not divisible by $p$, so $\nu(n) \in (A/\mathfrak{m})^{\times}$,
so $\nu(n)$ is a unit in $A$ since $A$ is local. So we can solve relation \eqref{hazewinkel relation 20} to get
\[ c_a = \frac{\gamma}{\nu(n)} (a^n - a)\]
for all $a\in A$. Hence $\overline{L}^A_{2n-2}$ is generated by $\gamma$.
\item {\bf If $n = p^t$ and either $A/\mathfrak{m}$ is infinite or $A/\mathfrak{m}$ has $p^s$ elements and $s\nmid t$:} 
Then there exists some element $a\in A/\mathfrak{m}$
such that $a^{p^t} \neq a$, and since $a$ is nonzero and $A$ is Henselian, $a$ lifts to an element
$\tilde{a}\in A$ such that $\tilde{a}^{p^t} - \tilde{a}$ is not in the maximal ideal in $A$.
Consequently $\tilde{a}^{p^t} - \tilde{a} \in A^{\times}$ and
hence we can solve relation \eqref{hazewinkel relation 20} to get
\[ \gamma = \frac{pc_{\tilde{a}}}{\tilde{a}^{p^t} - \tilde{a}}\]
and we can solve relation \eqref{hazewinkel relation 22} to get
\begin{equation*}\label{relation 176014} c_b = \frac{b^{p^t} - b}{\tilde{a}^{p^t} - \tilde{a}}c_{\tilde{a}},\end{equation*}
hence $\overline{L}^A_{2n-2}$ is generated by $c_{\tilde{a}}$. 
\item {\bf If $n = p^t$ and $A/\mathfrak{m}$ has $p^s$ elements and $s\mid t$:} 
This line of argument was inspired by Hazewinkel's Proposition 21.3.1 in \cite{MR506881}.
Let $M$ denote the $A$-submodule of $\overline{L}^A_{2n-2}$ generated by all the elements
$c_{m}$ with $m\in \mathfrak{m}$. 
Solving relation \eqref{hazewinkel relation 22}, we get
\begin{equation}\label{relation 176015} (a^{p^t}-a)c_m = (m^{p^t}-m)c_a \end{equation}
for all $a,m\in A$.
If $m\in \mathfrak{m}$, then $m^{p^t-1} - 1 \notin \mathfrak{m}$, hence $m^{p^t-1} - 1$ is a unit since $A$ is local.
Hence
\begin{equation}\label{relation 176016} \frac{a^{p^t} - a}{m^{p^t-1} - 1}c_m = mc_a \end{equation}
for all $m\in \mathfrak{m}$ and all $a\in A$ with $a\notin\mathfrak{m}$.
Relation \eqref{hazewinkel relation 20} yields
\begin{equation}\label{relation 176017} \frac{p}{1-m^{p^t-1}}c_m = m\gamma, \end{equation}
and now equation \eqref{relation 176016} and \eqref{relation 176017} give us that $\mathfrak{m}$ acts by zero on $\overline{L}^A_{2n-2}/M$, i.e.,
$\overline{L}^A_{2n-2}/M$ is an $A/\mathfrak{m}$-vector space.

Now $s$ divides $t$, and hence $x^{p^t} = x$ for all $x\in A/\mathfrak{m}$, so
relation \eqref{hazewinkel relation 21} becomes $c_{a+b} = c_a + c_b$ in $\overline{L}^A_{2n-2}/M$. Similarly, relation \eqref{hazewinkel relation 22}
becomes $c_{ab} = ac^b + bc^a$, i.e., the map 
\begin{align*}
  c: A/\mathfrak{m} &\rightarrow \overline{L}^A_{2n-2}/M \\\ 
  a&\mapsto c_a
\end{align*}
is a $\mathbb{Z}$-linear derivation. But the relevant module of K\"{a}hler differentials $\Omega^1_{(A/\mathfrak{m})/\mathbb{Z}}$ vanishes, since $A/\mathfrak{m}$ is a field. Hence $c$ factors through the zero module, i.e., $c$ is the zero map.
So $c_a = 0$ in $\overline{L}^A_{2n-2}/M$ for all $a\in A$ with $a\notin \mathfrak{m}$,
and $c_a = 0$ in $\overline{L}^A_{2n-2}/M$ for all $a\in  \mathfrak{m}$ by the definition of $M$.
So $\overline{L}^A_{2n-2}/M \cong A/\mathfrak{m}$ generated by $\gamma$. Consequently the elements $c_m$ for $m\in \mathfrak{m}$, along with $\gamma$, 
form a set of generators for
$\overline{L}^A_{2n-2}$.
\end{itemize}
\end{proof}

\begin{corollary}\label{hensel local ring lazard ring finite generation 2}
Let $A$ be a torsion-free Henselian local commutative ring with maximal ideal $\mathfrak{m}$, 
and suppose that $\mathfrak{m}$ is finitely generated.
Then, for each positive integer $n$,
$\overline{L}^A_{2n-2}$ is a finitely generated $A$-module.
\end{corollary}

\begin{corollary}\label{separatedness and completeness in complete case}
Let $A$ be a torsion-free Noetherian complete local commutative ring.
Then, for each integer $n$, the grading degree $n$ summand $L^A_n$
of the classifying ring $L^A$ of formal $A$-modules is
$\mathfrak{m}$-adically separated and $\mathfrak{m}$-adically
complete.
\end{corollary}
\begin{proof}
Krull's intersection theorem 
implies that every finitely generated $A$-module is $\mathfrak{m}$-adically separated, and it is elementary 
that every finitely generated
module over a Noetherian complete local ring with maximal ideal
$\mathfrak{m}$ is $\mathfrak{m}$-adically complete.
\end{proof}

Lemma \ref{coherent modules have exact completions} is immediate, when $R$ is Noetherian. The use of the lemma
is when $R$ is not Noetherian but $R^0$ is, e.g. $R\cong MU_* \cong L^{\mathbb{Z}}$.
\begin{lemma}\label{coherent modules have exact completions}
Let $R$ be a $\mathbb{Z}$-graded commutative ring which is
connective, i.e., there exists some integer $n$ such that $R^m\cong 0$
for all $m<n$. Assume furthermore that $R^0$ is Noetherian and that 
$R^i$ is a finitely generated $R^0$-module for each integer $i$. 

Then, for any $\mathbb{Z}$-graded finitely generated $R$-module $M$ 
and any ideal $I$ of $R$ generated by elements in $R^0$, 
the natural map 
\[\hat{R}_I\otimes_R M\rightarrow \hat{M}_I\]
is an isomorphism of $\mathbb{Z}$-graded $\hat{R}_I$-modules.
\end{lemma}
\begin{proof}
Since $M$ is finitely generated as an $R$-module, $M^i$ is finitely
generated as an $R^0$-module for any integer $i$, and since $M^i$ is
a finitely generated module over the Noetherian ring $R^0$, the map
$\hat{R}_I^0\otimes_{R^0} M^i\rightarrow \hat{M}^i_I$ is an isomorphism 
for all $i$.
\end{proof}

\begin{definition-proposition}\label{def of canonical ideal}
Let $A$ be a commutative ring and let $I$ be an ideal in $A$. Let $F$ be a formal group law defined over 
an $A$-algebra $R$. Equip the power series ring $R[[X]]$ with the $I+(X)$-adic filtration\footnote{To avoid any potential confusion, we remind the reader that the endomorphism ring $\End(F)$ of $F$ is a subset of $R[[X]]$, but the multiplication in $\End(F)$ is given by composition, and the addition in $\End(F)$ is given by formal addition using $F$. Consequently $\End(F)$ is not a subring of the power series ring $R[[X]]$.}, i.e., the decreasing filtration 
\begin{equation}\label{c-adic filt 0} R[[X]] = F_0 \supseteq F_1 \supseteq F_2 \supseteq \dots \end{equation}
in which $F_n$ is the $n$th power of the sum of the ideals $I$ and $(X)$. 
Let $\mathfrak{c}_n$ be the intersection of $F_n$ with $\End(F)\subseteq R[[X]]$. Then 
\begin{equation}\label{c-adic filt} \End(F) = \mathfrak{c}_0 \supseteq \mathfrak{c}_1 \supseteq \mathfrak{c}_2 \supseteq \dots \end{equation}
is a sequence of two-sided ideals of $\End(F)$. Furthermore, if $x\in \mathfrak{c}_m$ and $y\in \mathfrak{c}_n$, then $xy\in \mathfrak{c}_{m+n}$.

We refer to the filtration \eqref{c-adic filt} as the {\em $\mathfrak{c}$-adic filtration} on $\End(F)$. We call the resulting topology on $\End(F)$, in which the ideals \eqref{c-adic filt} are a neighborhood basis of zero,  the {\em $\mathfrak{c}$-adic topology}.
\end{definition-proposition}
\begin{proof}
Let $f(X),g(X)\in \End(F)$. 
It is routine to verify that
\begin{itemize} 
\item if $f(X),g(X)\in\mathfrak{c}_n$, then $f+g\in \mathfrak{c}_n$,
\item and if $f(X)\in \mathfrak{c}_m$ and $g(X)\in\mathfrak{c}_n$, then $fg\in \mathfrak{c}_{m+n}$ (note that this product in $\End(F)$ is the composition of power series).
\end{itemize}
The special cases $m=0$ and $n=0$ of the latter observation establish that each $\mathfrak{c}_n$ is a two-sided ideal of $\End(F)$.
\end{proof}

\begin{lemma}\label{separatedness and completeness lemma}
Let $A$ be a Noetherian commutative ring and let $I$ be an ideal in $A$. 
Let $R$ be a commutative $A$-algebra, and let $F$ be a formal $A$-module. 
If $R$ is $I$-adically separated, then $\End(F)$ is $\mathfrak{c}$-adically separated. If $R$ is $I$-adically complete, then $\End(F)$
is $\mathfrak{c}$-adically complete.
\end{lemma}
\begin{proof}
If $R$ is $I$-adically separated, then the filtration \eqref{c-adic filt 0} on $R[[X]]$ is as well, so the intersection 
\[ \bigcap_n \mathfrak{c}_n = \cap_n (F_n\cap \End(F))\]
must be zero.

Now suppose instead that $R$ is $I$-adically complete, 
and that 
\begin{equation}\label{cauchy seq 100}(\zeta_1(X), \zeta_2(X), \dots)\end{equation}
is a Cauchy sequence in the $\mathfrak{c}$-adic topology on $\End(F)$.
Let $\zeta_{m,n}$ denote the $n$th coefficient in the power series $\zeta_m(X)$. 
For each $n$, the sequence $\zeta_{1,n}, \zeta_{2,n}, \dots$ is an $I$-adically Cauchy sequence in $A$, hence converges. Hence the limit $\lim_m \zeta_m(X)$ exists in the $\mathfrak{c}$-adic topology: it is merely the endomorphism of $F$ whose $n$th power series coefficient is $\lim_m \zeta_{m,n}$.
So every $\mathfrak{c}$-adic Cauchy sequence in $\End(F)$ converges, so $\End(F)$ is $\mathfrak{c}$-adically complete.
\end{proof}

\begin{lemma}\label{extension of formal module structure via completion}
Let $A$ be a Noetherian commutative ring and let $I$ be an ideal in $A$. 
Suppose that $A$ is separated (but not necessarily complete) in the $I$-adic topology.
Let $R$ be a commutative $A$-algebra which is $I$-adically separated and complete, and let $F$ be a formal $A$-module over $R$.
Then $F$ is the underlying formal $A$-module of exactly one formal $\hat{A}_I$-module.
That is, the action map $\rho: A \rightarrow \End(F)$ extends uniquely to an action map
$\tilde{\rho}: \hat{A}_I\rightarrow \End(F)$ making $F$ a formal $\hat{A}_I$-module.
\end{lemma}
\begin{proof}
Choose an element $a\in \hat{A}_I$, and for each positive integer $n$, let $a_n$ be the image of $a$ under the
projection map $\hat{A}_I\rightarrow A/I^n$, and let $\tilde{a}_n$ be an element of $A$ whose reduction modulo $I^n$ is
$a_n$. (In other words: choose a sequence of elements $(\tilde{a}_1, \tilde{a}_2, \dots)$ of $A$ converging to 
$a$ in the $I$-adic topology.)
Then the sequence $(\tilde{a}_1, \tilde{a}_2, \dots )$ uniquely determines the element $a\in \hat{A}_I$,
since we assumed that $A$ is separated in the $I$-adic topology.

The tangent condition on $\rho$ (that $\rho(X) \equiv X \mod X^2$) implies that the image of $I^n$ under $\rho$ is contained in the ideal $\mathfrak{c}_n$ of Definition-Proposition \ref{def of canonical ideal}. 
By Lemma \ref{separatedness and completeness lemma}, $\End(F)$ is $\mathfrak{c}$-adically separated and complete,
so the sequence
$( \rho(\tilde{a}_1), \rho(\tilde{a}_2), \dots)$ in $\End(F)$ converges to a unique element in $\End(F)$. Let $\tilde{\rho}(a)$
be defined to be this element. It is elementary to check that the resulting map $\tilde{\rho}: \hat{A}_I \rightarrow \End(F)$
is a well-defined ring homomorphism and agrees with $\rho$ when composed with the injection $A\hookrightarrow \hat{A}_I$.
(This map $\tilde{\rho}$ is, of course, the one given by the universal property of completion, but we are giving some detail here because
$\End(F)$ is not typically commutative and $\rho$ does not typically have its image inside the center of $\End(F)$, so the situation
is not exactly the textbook one encountered in algebra.)
The tangent condition for $\tilde{\rho}$ is similarly easy: any element $a\in \hat{A}_I$ can be approximated arbitrarily 
$\mathfrak{c}$-adically closely by an element of $A$, and $\tilde{\rho}$ satisfies the tangent condition on elements of $A$
since $\tilde{\rho}$ coincides with $\rho$ on elements of $A$. 

Consequently $F$ is indeed the underlying formal $A$-module of a formal $\hat{A}_I$-module.
The fact that the ring homomorphism $\tilde{\rho}$ is the {\em unique} extension of $\rho$ to a ring map $\hat{A}_I\rightarrow \End(F)$ 
is as follows: any ring homomorphism $\hat{A}_I\rightarrow \End(F)$ extending $\rho$ sends each $I^n$ into $\mathfrak{c}_n$ and hence is continuous,
hence is a continuous homomorphism of abelian groups;
now the universal property of the completion implies that the extension $\tilde{\rho}$ is unique.
\end{proof}

\begin{lemma}\label{strict isos under completion}
Let $A$ be a Noetherian commutative ring and let $I$ be an ideal in $A$.
Suppose that $A$ is $I$-adically separated, and suppose that $R$ is an $I$-adically separated and complete $A$-algebra. Let $F,G$ be formal $\hat{A}_I$-modules over $R$.
Suppose that $f: F \rightarrow G$ is a strict isomorphism of the underlying formal $A$-modules of $F$ and $G$.
Then $f$ is also a strict isomorphism $F\rightarrow G$ of formal $\hat{A}_I$-modules.
\end{lemma}
\begin{proof}
Let $\rho_F: \hat{A}_I\rightarrow \End(F)$ and $\rho_G:\hat{A}_I \rightarrow \End(G)$ denote the structure maps of
$F$ and $G$ as formal $\hat{A}_I$-modules, respectively. (The existence and uniqueness of these structure maps follows from Lemma \ref{extension of formal module structure via completion}.) Then $\rho_F(f(a)(X)) = f(\rho_G(a)(X))$ for all $a\in A$, and we need
to show that the same is true for all $a\in \hat{A}_I$. For any $a\in \hat{A}_I$, choose a sequence of elements $a_1, a_2, \dots $ in $A$ such
that $\lim_{n\rightarrow\infty} a_n = a$ in the $I$-adic topology. Then the fact that $\rho_F$ and $\rho_G$ are continuous (since each sends
$I$ into $\mathfrak{c}$) implies that
\begin{align}
\label{equality one billion} f(\rho_F(\lim_{n\rightarrow \infty} a_n)(X)) &= 
  f(\lim_{n\rightarrow\infty} \rho_F(a_n)(X)) \mbox{,\ \ and} \\
\label{equality one billion and one} \rho_G(\lim_{n\rightarrow \infty} a_n)(f(X)) &=
  \lim_{n\rightarrow\infty} \rho_G(a_n)(f(X)).
\end{align}

Now since $f$ is a homomorphism of formal group laws, its constant coefficient is zero, 
and hence $f$ is continuous in a limited sense:
whenever $\xi_1(X), \xi_2(X), \dots $ is a $\mathfrak{c}$-adically convergent 
sequence in $\End(F)$ such that each power series \linebreak $f(\xi_1(X)), f(\xi_2(X)), \dots$ is contained
in $\End(G)\subseteq R[[X]]$ and $\mathfrak{c}$-adically convergent,
we get an equality $\lim_{n\rightarrow\infty} f(\xi_n(X)) = f(\lim_{n\rightarrow\infty} \xi_n(X))$.
Consequently \eqref{equality one billion} is equal to \eqref{equality one billion and one},
and hence $\rho_F(f(a)(X)) = f(\rho_G(a)(X))$.
\end{proof}

\begin{theorem}\label{completion of the hopf algebroid}
Let $A$ be a torsion-free commutative ring which is finitely generated as a commutative ring. 
Let $I$ be a maximal ideal of $A$, and suppose that $A$ is separated\footnote{This separation condition is automatically satisfied if $A$ is an integral domain, by Krull's intersection theorem.} in the $I$-adic topology.
Then the natural maps of graded Hopf algebroids
\begin{align} 
\label{completion map 344} (L^A\otimes_{A} \hat{A}_I , L^AB \otimes_{A} \hat{A}_I ) 
  &\rightarrow ((L^A)^{\hat{}}_I , (L^AB)^{\hat{}}_I ) \\
\label{completion map 345}  &\rightarrow (L^{\hat{A}_I} , L^{\hat{A}_I}B ) \end{align}
are isomorphisms.
\end{theorem}
\begin{proof}
By Proposition \ref{finite typeness of lazard ring}, for all integers $n$ the degree $n$ summand $L^A_n$ in the ring $L^A$ 
is a finitely-generated $A$-module, and 
$L^AB \cong L^AB\otimes_{L} LB\cong L^A[t_1, t_2, \dots ]$ 
(by Proposition \ref{formal module base change thm} and Theorem \ref{basic existence thm on classifying hopf algebroids}) is also a finitely-generated $A$-module in each degree. 
Consequently Lemma \ref{coherent modules have exact completions} applies, since $A$ is finitely generated as a commutative ring and hence Noetherian, even though $L^A$ is typically not Noetherian. So the map \eqref{completion map 344} is an isomorphism.

The more substantial result is that \eqref{completion map 345} is also an isomorphism.
By Lemma \ref{extension of formal module structure via completion}, Lemma \ref{strict isos under completion}, and the
universal properties of the rings involved, 
the map \eqref{completion map 345} induces bijections
\begin{align*}
 \hom_{\hat{A}_I-alg}(L^{\hat{A}_I}, R) &\stackrel{\cong}{\longrightarrow} \hom_{\hat{A}_I-alg}(L^{A}\otimes_A \hat{A}_I, R)\mbox{\ \ and} \\
 \hom_{\hat{A}_I-alg}(L^{\hat{A}_I}B, R) &\stackrel{\cong}{\longrightarrow} \hom_{\hat{A}_I-alg}(L^{A}B \otimes_A \hat{A}_I , R),
\end{align*}
natural in $R$, for all commutative $\hat{A}_I$-algebras $R$ which are $I$-adically separated and complete.

Now the Yoneda lemma tells us that the ring maps
$L^A \otimes_A \hat{A}_I \rightarrow L^{\hat{A}_I}$ and $L^AB \otimes_A \hat{A}_I \rightarrow L^{\hat{A}_I}B$
are isomorphisms, as long as all four of these rings are actually objects in the category of
$I$-adically separated and complete commutative $\hat{A}_I$-algebras!
The graded ring $L^{A_I}$ a finitely generated $A$-module in each degree by Corollary \ref{finite generation and separatedness of local lazard ring}. Hence the same is true of $(L^{A_I})^{\hat{}}_I \cong L^{A_I}\otimes_{A_I}\hat{A}_I$, hence
$(L^{A_I})^{\hat{}}_I$ is $I$-adically separated and $I$-adically complete
by the same argument as in the proof of Corollary \ref{separatedness and completeness in complete case},
and the same is true for $(L^{A_I}B)^{\hat{}}_I \cong (L^{A_I})^{\hat{}}_I \otimes_L LB$, by Proposition \ref{formal module base change thm}.
On the other hand, $L^{\hat{A}_I}$ is $I$-adically separated and complete in each grading degree by Corollary \ref{separatedness and completeness in complete case}, hence $L^{\hat{A}_I}B \cong L^{\hat{A}_I}\otimes_L LB$ 
is as well, again by Proposition \ref{formal module base change thm}.
\end{proof}

\begin{corollary}\label{completion and cohomology}
Let $A$ be a torsion-free commutative ring which is finitely generated as a commutative ring.
Let $I$ be a maximal ideal of $A$, and suppose that $A$ is separated in the $I$-adic topology. 
Let $M$ be a graded left $L^{A}B$-comodule which is finitely-generated
as an $A$-module in each degree,
and suppose that $M$ is bounded-below, i.e., there exists some integer $b$
such that $M$ is trivial below degree $b$.
(For example, $M = L^A$ satisfies all these conditions on $M$.)
Then, for all integers $s,t$ with $s\geq 0$, we have 
isomorphisms of $\hat{A}_I$-modules
\begin{align} 
\Ext_{(L^A, L^AB)}^{s,t}(L^A, M) \otimes_A \hat{A}_I
 \label{iso 1892}  &\cong \Ext_{(L^A, L^AB)}^{s,t}(L^A, \hat{M}_I) \\
 \label{iso 1894}  &\cong \Ext_{(L^{\hat{A}_I}, L^{\hat{A}_I}B)}^{s,t}(L^{\hat{A}_I}, \hat{M}_I) \end{align}
\end{corollary}
\begin{proof}
Let $C^{\bullet}_{(L^A, L^AB)}(M)$ be the cobar complex\footnote{See Appendix 1 of \cite{MR860042} for the definition and basic properties of the cobar complex. All that is necessary for us to know right now about the cobar complex is that its module of $n$-cochains $C^n_{(L^A,L^AB)}(M)$ is isomorphic to $(L^AB)^{\otimes_{L^A} n}\otimes_{L^A} M$, and the cohomology of $C^{\bullet}_{(L^A,L^AB)}(M)$ is $\Ext^{*,*}_{(L^A,L^AB)}(L^A, M)$, i.e., $\Cotor^{*,*}_{L^AB}(L^A,M)$.} of the Hopf algebroid $(L^A, L^AB)$ with coefficients in $M$.
Then, since $L^AB$ is a finitely generated $A$-module in 
each degree by Proposition \ref{finite typeness of lazard ring},
the same is true of $L^AB\otimes_{L^A}L^AB \otimes_{L^A}\dots \otimes_{L^A} M = (L^AB)^{\otimes_{L^A} n}\otimes_{L^A} M$.
Consequently we have isomorphisms
\begin{align*} 
 (L^AB)^{\otimes_{L^A} n}\otimes_{L^A} \hat{M}_I  
  &\cong (L^AB)^{\otimes_{L^A} n}\otimes_{L^A} M\otimes_A \hat{A}_I \\
  &\cong \left(L^AB\otimes_A \hat{A}_I\right)^{\otimes_{L^A\otimes_A \hat{A}_I} n}\otimes_{L^A\otimes_A \hat{A}_I} \left(M \otimes_A \hat{A}_I\right) \\
  &\cong (L^{\hat{A}_I}B)^{\otimes_{L^{\hat{A}_I}} n}\otimes_{L^{\hat{A}_I}} \hat{M}_I,\end{align*}
which is the module of $n$-cochains $C^n_{(L^{\hat{A}_I}, L^{\hat{A}_I}B)}(\hat{M}_I)$.
These isomorphisms are natural, commuting with the cobar complex differentials,
hence giving us isomorphism \eqref{iso 1894}.
Meanwhile, isomorphism \eqref{iso 1892} follows from $\hat{A}_I$ being a flat
$A$-module (classical, as in Proposition 10.14 in \cite{MR0242802}),
hence tensoring with $\hat{A}_I$ commutes with taking cohomology of the cobar complexes.
\end{proof}

\begin{corollary}\label{ss completion corollary}
Let $A,I,M$ be as in Corollary \ref{completion and cohomology}. Let $E_0^IM$ denote the associated graded comodule of the $I$-adic filtration on $M$.
Then there exists a conditionally convergent spectral sequence
\begin{align}
\nonumber E_1^{s,t,u} \cong \Cotor_{L^AB}^{s,t,u}\left( L^A, E_0^IM\right)
 &\Rightarrow \Cotor^{s,t}_{L^AB}\left( L^A,M \right)^{\hat{}}_I \\
\label{iso 3234094} &\cong \Cotor^{s,t}_{L^{\hat{A}_I}B}\left( L^{\hat{A}_I},\hat{M}_I \right) \\
\nonumber d_r: E_r^{s,t,u} &\rightarrow E_r^{s+1,t,u+r}.
\end{align}
\end{corollary}
\begin{proof}
This is the spectral sequence of the $I$-adic filtration on the cobar complex of $(L^A,L^AB)$ with coefficients in $M$. The isomorphism \eqref{iso 3234094} is due to Corollary \ref{completion and cohomology}.
\end{proof}
See \cref{Construction of the 2-adic...} and \cref{cotor for Z} for examples of explicit calculations with the spectral sequence of Corollary \ref{ss completion corollary}. 

\subsection{Point-detecting ideals and the fundamental functional.}
\label{Point-detecting ideals...}

\begin{definition}
Let $A$ be a commutative ring, and let $k$ be a field. We say that an ideal $I$ of $A$ {\em detects all $k$-points} $I$ is in the kernel of every ring homomorphism $A \rightarrow k$.
\end{definition}

Let $k$ be the finite field with $p^n$ elements. Then there exists a {\em largest} ideal of $A$ which detects all $k$-points: namely, the ideal generated by $p$ and by the difference $a^{p^n} - a$ for each $a\in A$. We will write $I_{p^n}$ for the largest ideal of $A$ which detects all $k$-points, so that
\[ I_{p^n} = \left( p, a^{p^n} - a\ \ \mbox{for\ all}\ a\in A\right).\]
We call $I_{p^n}$ the {\em universal $\mathbb{F}_{p^n}$-point-detecting ideal} of $A$.

\begin{prop}\label{honest point-counting}
Let $A$ be a Noetherian commutative ring. Then the quotient map $A\rightarrow A/I_p$ coincides with the universal map $A\rightarrow \prod \mathbb{F}_p$, with the product taken over all ring homomorphisms $A\rightarrow \mathbb{F}_p$. Put another way: $\Spec A/I_p$ is the union of the $\mathbb{F}_p$-points of $\Spec A$. 
\end{prop}
\begin{proof}
Elementary, but here is the argument: clearly $A/I_p$ has no nonzero nilpotents, so $I_p$ is reduced. Consequently, in its primary decomposition $I_p = \mathfrak{p}_1 \cap \dots \cap \mathfrak{p}_n$, each primary ideal $\mathfrak{p}_j$ is prime. The quotient map $A\rightarrow A/\mathfrak{p}_j$ must factor through $A\rightarrow A/I_p$, so $A/\mathfrak{p}_j$ is an integral domain which is a quotient of $A/I_p$. In such an integral domain, we have $a(a^{p-1} - 1) = 0$ for all $a$, hence $a^{p-2} = a^{-1}$ for all nonzero $a$, i.e., $A/\mathfrak{p}_j$ is a field. Hence $I_p$ is an intersection of finitely many coprime maximal ideals, so the Chinese Remainder Theorem ensures that $A/I_p$ is a product of fields.
\end{proof}

\begin{remark} \label{remark on local zeta}
The local zeta-function of an affine variety $V = \Spec A$ at some prime number $p$ is defined as
\[ Z(V,t) = exp\left(\sum_{m\geq 1} \frac{N_m}{m} t^m\right),\]
where $N_m$ is the number of $\mathbb{F}_{p^m}$-points of $V$. Proposition \ref{honest point-counting} shows that the order of $A/I_p$ is $p^{N_1}$. An analogous statement holds for other prime powers, with appropriate adjustments: $A/I_{p^2}$ ``counts'' both the $\mathbb{F}_p$-points and the $\mathbb{F}_{p^2}$-points of $\Spec A$, for example.
\end{remark}

\begin{definition}\label{def of fundamental functional}
Let $n$ be a positive integer, and let $A$ be a commutative ring which is $\nu(n)$-torsion-free. Recall from Proposition \ref{drinfeld presentation} 
that $\overline{L}^A_{2n-2}$ is described by Drinfeld's presentation:
it is generated, as an $A$-module, by elements $\gamma$ and $\{ c_a\}_{a\in A}$, subject to the relations \eqref{hazewinkel relation 20}, \eqref{hazewinkel relation 21}, and \eqref{hazewinkel relation 22}.

Let $M^A_{2n-2}$ denote the $A$-module generated by elements $\gamma$ and $\{ c_a\}_{a\in A}$, subject only to the relation \eqref{hazewinkel relation 20}, i.e.,
\begin{align*}
 (a^n - a)\gamma &= \nu(n)c_a \mbox{\ \ for\ all\ } a\in A,\end{align*}
where $\nu(n)$ is as in Proposition \ref{drinfeld presentation}. That is, $\nu(n) = p$ if $n$ is a power of a prime $p$, and $\nu(n) = 1$ if $n$ is not a prime power.
Let $q^A_{2n-2}: M^A_{2n-2} \rightarrow \overline{L}^A_{2n-2}$ denote the obvious $A$-module quotient map.

By the {\em degree $2n-2$ fundamental functional of $A$}, we mean the unique $A$-module homomorphism $\sigma_{2n-2}: \overline{L}^A_{2n-2} \rightarrow A$ given by
\begin{align*} 
 \sigma_{2n-2}(\gamma) &= \nu(n),\mbox{ \ \ and}\\
 \sigma_{2n-2}(c_a) &= a^n-a.
\end{align*}
\end{definition}

\begin{prop}\label{point-detecting iso}
Let $A$ be a commutative ring. Suppose that $A$ has no nonzero $2$-torsion.
Then the degree $2$ fundamental functional $\sigma_2: \overline{L}^A_2\rightarrow A$ is injective, and its image is the universal $\mathbb{F}_{2}$-point-detecting ideal $I_2$ of $A$. Consequently we have an isomorphism of $A$-modules $\overline{L}^A_2 \stackrel{\cong}{\rightarrow} I_2$.
\end{prop}
\begin{proof}
It is straightforward from the definition of $\sigma_2$ that its image is generated by $2$ and by $a^2-a$ for each $a\in A$, i.e., $I_2 = \im\sigma_2$. All that needs to be checked is that $\sigma_2$ is injective. 

Our argument for injectivity of $\sigma_2$ involves the Hochschild homology of the commutative ring $A/2$, with coefficients in a particular $A/2$-bimodule $\Psi$. To be clear, by $A/2$ we mean the quotient of the ring $A$ by its principal ideal generated by $2$. The $A/2$-bimodule $\Psi$ is defined as follows: 
\begin{itemize}
\item as a left $A/2$-module, $\Psi$ is simply $A/2$ itself, i.e., free on one generator,
\item and the right $A/2$-action on $\Psi$ is given by letting $a\cdot x$ be the product $a^2x$ in the ring $A/2$.
\end{itemize}
The cyclic bar complex (i.e., the standard resolution for calculating Hochschild homology) with coefficients in $\Psi$ is as follows:
\[
 0 \leftarrow
  A/2 \stackrel{d_0}{\longleftarrow} 
 A/2\otimes_{\mathbb{F}_2} A/2 \stackrel{d_1}{\longleftarrow}
  A/2\otimes_{\mathbb{F}_2} A/2\otimes_{\mathbb{F}_2} A/2 \stackrel{d_2}{\longleftarrow}
  \dots
\]
with
\begin{align*}
 d_0( a_0\otimes a_1 ) &= (a_1^2 - a_1)a_0 \mbox{,\ \ and} \\
 d_1( a_0\otimes a_1\otimes a_2 ) &= a_0a_1^2\otimes a_2 - a_0 \otimes a_1a_2 + a_2a_0\otimes a_1 .
\end{align*}
Let $P^A$ be the cokernel of the $A$-module homomorphism $A\rightarrow \overline{L}^A_2$ sending $1$ to $\gamma$.
Consider the kernel $\ker d_0$ as a left $A/2$-module. It admits an $A/2$-module homomorphism $f: \ker d_0 \rightarrow P^A$ given by $f(a_0\otimes a_1) = a_0c_{a_1}$.
Let $\theta: P^A \rightarrow A/2$ be the $A/2$-module map given by $\theta(c_a) = a^2-a$.
It is routine to check that the resulting sequence of $A/2$-modules
\[ 0 \rightarrow \im d_1 \rightarrow \ker d_0 \stackrel{f}{\longrightarrow} P^A \stackrel{\theta}{\longrightarrow} A/2 \] 
is exact, by checking that Drinfeld relation \eqref{hazewinkel relation 22} is precisely the relation imposed on $\ker d_0$ by quotienting out by its submodule $\im d_1$. 
Consequently the Hochschild homology group $HH_1(A/2; \Psi)$ vanishes if and only if $\theta$ is injective.

Now consider the commutative diagram with exact rows
\begin{equation}\label{comm diag 09593}\xymatrix{
 \ker \gamma \ar[r]\ar[d] & A \ar[d]^{\id} \ar[r]^(.4){\gamma} & \overline{L}^A_{2}\ar[d]^{\sigma_2}\ar[r] & P^A\ar[r] \ar[d]^{\theta} & 0 \ar[d] \\
 0 \ar[r] & A \ar[r]^{2} & A \ar[r] & A/2 \ar[r] & 0 .}\end{equation}
A routine diagram chase in diagram \eqref{comm diag 09593} shows that $\sigma_2$ is injective if and only if $\theta$ is injective.

Finally, we invoke a 2007 calculation of Pirashvili, the main result of \cite{MR2336249}: if $n$ is a positive power of a prime number $p$, $B$ is a commutative $\mathbb{F}_p$-algebra, and $\Phi^n(B)$ is the ring $B$ regarded as a $B$-bimodule via the free $B$-action on the left and via the Frobenius-twisted $B$-action $x\cdot b := b^p x$ on the right, then the Hochschild homology groups $HH_i(B; \Phi^n(B))$ are trivial for all $i>0$. In particular, our Hochschild group $HH_1(A/2; \Phi)$ is trivial. Hence $\theta$ is injective, hence $\sigma_2$ is injective, as desired.
\end{proof}
Proposition \ref{point-detecting iso} admits a generalization to the higher fundamental functionals $\sigma_n$ and the universal $k$-point-detecting ideals for all finite fields $k$. That generalization is not used in this paper, so we do not present it here, but it can be found in the preprint \cite{cmah2}.

\section{Generalities on the moduli of formal $A$-module $2$-buds.}
\label{Generalities on the moduli...}

Now we begin to restrict our attention from formal $A$-modules down to only their $2$-buds, in order to facilitate explicit calculations which hold for a very large class of rings $A$.

\subsection{Structure of the classifying ring of formal $A$-module $2$-buds.}
\label{Moduli...2-buds}

A formal group law $2$-bud is given by a single coefficient $\gamma$, so that $F(X,Y) \equiv X+Y + \gamma XY\mod (X,Y)^3$. The $A$-action map $\rho: A \rightarrow \End(F)$ is given by $\rho_a(X) \equiv aX + c_{a}X^2\mod X^3$. It is classical (see \cite{MR0384707}) that, by elementary calculation, one finds that the relations among the elements $\gamma$ and $c_{a}$ imposed by associativity of $F$, by $\rho_{ab}(X) = \rho_a(\rho_b(X))$, by $\rho_{a+b}(X) = F(\rho_a(X),\rho_b(X))$, and by $F(\rho_a(X),\rho_a(Y)) = \rho_aF(X,Y)$, are the Drinfeld relations
\begin{align}
\label{drinfeld rel 1} (a^2-a)\gamma &= 2c_{a} \\
\label{drinfeld rel 2} c_{a+b} - c_{a} - c_{b} &= ab\gamma \\
\label{drinfeld rel 3} ac_{b} + b^2c_{a} &= c_{ab}.
\end{align}
\begin{definition} Let $\Dr^A$ denote the $A$-module generated by the symbols $\gamma$ and $c_{a}$ for each $a\in A$, subject to the relations \eqref{drinfeld rel 1} through \eqref{drinfeld rel 3}. 
\end{definition} That is, $\Dr^A \cong \overline{L}^A_2$. It follows easily\footnote{The same claim cannot safely be made for formal module $n$-buds for $n>2$. Consider the situation in the case of $3$-buds: we have a parameter $\gamma_2$, and a parameter $c_{a,2}$ for each $a\in A$, which jointly determine the quadratic terms in the $3$-bud. We also have a deformation parameter $\gamma_3$, and a deformation parameter $c_{a,3}$ for each $a\in A$, for the cubic terms in the $3$-bud. The $A$-module $\overline{L}^A_3$ generated by $\gamma_3$ and $c_{a,3}$, related by the cubic Drinfeld relations, is not necessarily (i.e., not for all $A$) a {\em projective} $A$-module. What happens as a consequence is that the classifying ring $L^A_{3-buds}$ of formal $A$-module $3$-buds is not necessarily a symmetric $A$-algebra on $\overline{L}^A_2\oplus \overline{L}^A_3$. We have made some calculations of the ring $L^A_{3-buds}$ as a function of $A$, but those calculations do not fit within the scope of this paper.} that the classifying ring $L^A_{2-buds}$ of formal $A$-module $2$-buds is the symmetric $A$-algebra $\Symm_A(\Dr^A)$.

\begin{remark}
Among the elementary consequences of \eqref{drinfeld rel 2} and \eqref{drinfeld rel 3}, we have that $c_0 = 0 = c_1$, and that $c_{-1} = \gamma$. Consequently a smaller presentation for $\Dr^A$ is possible: $\Dr^A$ is the $A$-module with a generator $c_a$ for each $a\in A$, subject to the relations \eqref{drinfeld rel 1} through \eqref{drinfeld rel 3} but with $c_{-1}$ written in place of $\gamma$ throughout. This is convenient for the sake of understanding the moduli of formal $A$-module $2$-buds, but the coincidence $c_{-1} = \gamma$ is specific to the case of $2$-buds. Consequently it is easier to recognize patterns that hold for formal $A$-modules $n$-buds for $n>2$ (and for formal $A$-modules {\em simpliciter}) if we use the non-minimal presentation for $\Dr^A$, rather than replacing $\gamma$ with $c_{-1}$.
\end{remark}

\subsection{Structure of the moduli stack of formal $A$-module $2$-buds.}
\label{Moduli...2-buds 2}

It follows from \cref{Moduli...2-buds} and from the structure theory in \cref{generalities on LA and LAB} that the classifying Hopf algebroid $(L^A_{2-buds},L^A_{2-buds}B)$ of formal $A$-module $2$-buds has the following simple form:
\begin{align}
\nonumber L^A_{2-buds} &\cong \Symm_A(\Dr^A) ,\\
\label{iso 445010} L^A_{2-buds}B &\cong \Symm_A(\Dr^A)[t] ,\\
\nonumber \eta_L(x) &= x\ \ \mbox{for\ all\ } x\in L^A_{2-buds},\\
\label{iso 445011} \eta_R(\gamma_{}) &= \gamma_{} + 2t,\\
\label{iso 445012} \eta_R(c_{a}) &= c_{a} + (a^2-a)t \ \ \mbox{for all\ } a\in A,\\
\label{iso 445013} \Delta(t) &= t\otimes 1 + 1\otimes t,\\
\label{iso 445014} \epsilon(t) &= 0.
\end{align}
For any formal $A$-module $2$-bud $F$, the truncated power series $f(X) \equiv X + tX^2$ modulo $X^3$ is a strict isomorphism from $F$ to {\em some} formal $A$-module $2$-bud $G$. The parameter $t$ in \eqref{iso 445010} is the parameter $t$ in the strict isomorphism $f(X) = X + t X^2$ of formal $A$-module $2$-buds. If $F(X,Y) \equiv X+Y+\gamma_{}XY$ modulo $(X,Y)^3$, then $f$ is a strict isomorphism from $G$ to $G(X,Y) \equiv X+Y+(\gamma_{} + 2t)XY$ modulo $(X,Y)^3$. This yields \eqref{iso 445011}. 

Similarly, if the action map $\rho_F: A \rightarrow \End(F)$ of the formal $A$-module $2$-bud $F$ is given by $\rho_F(a)(X) = aX + c_{a}X^2$ modulo $X^3$, then the codomain of the strict isomorphism $f$ with domain $F$ has $A$-action map $\rho_G(a)(X) = aX + (c_{a} + (a^2-a)t)X^2$ modulo $X^3$. This yields \eqref{iso 445012}.

The formulas \eqref{iso 445013} and \eqref{iso 445014} arise from the observation that the composite of strict isomorphisms $f(X) \equiv X + aX^2$ modulo $X^3$ and $g(X) \equiv X + a^{\prime}X^2$ modulo $X^3$ is $(f\circ g)(X) \equiv X + (a+a^{\prime})X^2$ modulo $X^3$. 

Smoothness of the ring map $\eta_L$ implies that the stackification of the groupoid scheme represented by $(\Spec L^A_{2-buds},\Spec L^A_{2-buds}B)$ is an Artin stack. We write $\mathcal{M}_{fmA}^{2-buds}$ for the resulting {\em moduli stack of formal $A$-module $2$-buds}. To be clear, the groupoid of $R$-points of $\mathcal{M}_{fmA}^{2-buds}$ is the groupoid of $A$-module structures, with identity element $0$, on the first-order neighborhood of $\Spec R$ in $\mathbb{A}^1_R$. Such an $A$-module structure is slightly less data than a formal $A$-module $2$-bud. To go from the $A$-module structure to a formal $A$-module $2$-bud, one needs to specify an isomorphism of $\Spec R[X]/X^2$ with the first-order neighborhood of $\Spec R$ in $\mathbb{A}^1_R$. The situation is entirely analogous to that described in \cref{Formal modules and formal module n-buds.}.

As described in \cref{Hopf algebroids and stacks.}, it follows from standard arguments that the flat cohomology group $H^s_{fl}(\mathcal{M}_{fmA}^{2-buds}; \omega^{\otimes t})$ is isomorphic to $\Cotor^{s,2t}_{L^A_{2-buds}B}(L^A_{2-buds},L^A_{2-buds})$. Here $\omega$ is the line bundle of differentials, which coincides with the quasicoherent $\mathcal{O}_{\mathcal{M}_{fmA}^{2-buds}}$-module corresponding to the graded $L^AB$-comodule $\Sigma^2 L^A$. Since $L^A_{2-buds}$ and $L^A_{2-buds}B$ are trivial in odd degrees, the Cotor-groups \linebreak $\Cotor^{s,t}_{L^A_{2-buds}B}(L^A_{2-buds},L^A_{2-buds})$ are trivial in odd internal degrees $t$. Hence each of the nontrivial Cotor-groups $\Cotor^{s,t}_{L^A_{2-buds}B}(L^A_{2-buds},L^A_{2-buds})$ agrees with a flat cohomology group of $\mathcal{M}_{fmA}^{2-buds}$. Consequently, throughout the rest of this paper, we will often state our results in terms of flat cohomology of the Artin stack $\mathcal{M}_{fmA}^{2-buds}$, but the proofs will be carried out in terms of $\Cotor$ over the Hopf algebroid $(L^A_{2-buds}, L^A_{2-buds}B)$. The reader who strongly prefers one perspective over the other can translate all statements about flat cohomology into statements about $\Cotor$, or conversely.

\begin{lemma}\label{rel inj lemma 12}
Suppose $2$ is a unit in the commutative ring $A$. Then the $L^A_{2-buds}B$-comodule $L^A_{2-buds}$ is a summand of the $L^A_{2-buds}B$-comodule $L^A_{2-buds}B$.
\end{lemma}
\begin{proof}
By the assumption $\frac{1}{2}\in A$, the right unit map 
\begin{align*}
 \eta_R: L^A_{2-buds} & \rightarrow L^A_{2-buds}B = L^A_{2-buds}[t]   
\end{align*}
admits the retraction given by the $L^A_{2-buds}$-algebra map
\begin{align*}
 r: L^A_{2-buds}B & \rightarrow L^A_{2-buds} \\
 t &\mapsto \frac{d}{2}.
\end{align*}
It is elementary to verify that $r$ is a left $L^A_{2-buds}B$-comodule algebra homomorphism, and that $r\circ\eta_R = \id_{L^A_{2-buds}}$, i.e., $L^A_{2-buds}$ is a $L^A_{2-buds}B$-comodule retract of $L^A_{2-buds}B$.
\end{proof}

Here is a simple cohomological application of Theorem \ref{lazard ring localization iso}: 
\begin{prop}\label{2-loc prop}
Let $A$ be a commutative ring. Then, for all $s>0$ and all $t$, the abelian group 
$H^s_{fl}(\mathcal{M}_{fmA}^{2-buds};\omega^{\otimes t})$ is $2$-power-torsion. 
\end{prop}
\begin{proof}
By Theorem \ref{lazard ring localization iso}, we have
\begin{align}
\label{cotor iso 2843} \Cotor^{s,t}_{L^A_{2-buds}B}(L^A_{2-buds},L^A_{2-buds})[2^{-1}]
  &=  \Cotor^{s,t}_{L^{A[2^{-1}]}_{2-buds}B}(L^{A[2^{-1}]}_{2-buds},L^{A[2^{-1}]}_{2-buds}).
\end{align}
By Lemma \ref{rel inj lemma 12}, $L^{A[2^{-1}]}_{2-buds}$ is a relatively injective $L^{A[2^{-1}]}_{2-buds}B$-comodule (see appendix 1 of \cite{MR860042} for basic properties of relatively injective comodules). Consequently the right-hand side of \eqref{cotor iso 2843} is trivial for $s>0$.
Hence $\Cotor^{s,t}_{L^A_{2-buds}B}(L^A_{2-buds},L^A_{2-buds})$ is killed by inverting $2$ for $s>0$, i.e., $\Cotor^{s,t}_{L^A_{2-buds}B}(L^A_{2-buds},L^A_{2-buds})$ is $2$-power-torsion for $s>0$.
\end{proof}

\subsection{$H^*_{fl}(\mathcal{M}_{fmA}^{2-buds};\omega^{\otimes *})$ as the cohomology of a $\mathbb{G}_a$-action.}
\label{...as the cohomology of a...}

Here is a slightly more conceptual way to phrase the presentation for the Hopf algebroid given by \eqref{iso 445010} through \eqref{iso 445014}. Given a Hopf algebroid $(B,\Gamma)$ and a left $\Gamma$-comodule algebra $C$, one can form the Hopf algebroid $(C, \Gamma\otimes_B C)$ by the construction given in Proposition \ref{tensoring on one side}. We will refer to $(C, \Gamma\otimes_B C)$ as the {\em one-sided base-change Hopf algebroid} of $(B,\Gamma)$ along $B\rightarrow C$.
The left $\Gamma$-coaction on $C$ is necessary in order to define the right unit map on the one-sided base change Hopf algebroid, and making different choices of left $\Gamma$-coaction on $C$ will yield different Hopf algebroids $(C, \Gamma\otimes_B C)$. Specializing to the case of formal $A$-module $2$-buds: consider the situation where the Hopf algebroid $(B,\Gamma)$ is the Hopf $\mathbb{Z}$-algebr{\em a} representing the additive group scheme $\mathbb{G}_a$ over $A$. Concretely, $(B,\Gamma)$ is the Hopf algebra $(A,A[t])$, with $t$ primitive. We have a left $\Gamma$-coaction on the ring $L^A_{2-buds}$ given by the $A$-algebra map
\begin{align}
\label{coaction map 1} 
\psi: L^A_{2-buds} &\rightarrow A[t]\otimes_{A} L^A_{2-buds} \\
\nonumber d &\mapsto 1 \otimes d + 2t \otimes 1\\
\nonumber c_a &\mapsto 1 \otimes c_a + (a^2-a)t\otimes 1. 
\end{align}
The map $\psi$ is simply the ring map classifying the underlying formal $A$-module of the {\em target} of the universal strict isomorphism of formal $A$-modules. From \eqref{iso 445010} through \eqref{iso 445014}, we see that $(L^A_{2-buds}, L^A_{2-buds}B)$ is simply the one-sided base change of $(A,A[t])$ along the map $\psi$. 

One useful consequence, which we use in computations throughout the rest of this paper, is an isomorphism in $\Cotor$: by the change-of-rings isomorphism of Proposition \ref{tensoring on one side}, the identification of $(L^A_{2-buds}, L^A_{2-buds}B)$ as a one-sided base change Hopf algebroid yields isomorphisms of bigraded $A$-modules
\begin{align}
\label{change-of-rings iso 230}
 H^s_{fl}(\mathcal{M}_{fmA}^{2-buds};\omega^{\otimes t}) &\cong 
 \Cotor_{L^A_{2-buds}B}^{s,2t}(L^A_{2-buds},L^A_{2-buds})\\
  &\cong \Cotor_{A[t]}^{s,2t}(A,L^A_{2-buds}).
\end{align}

\subsection{Calculation of $H^0_{fl}(\mathcal{M}_{fmA}^{2-buds};\omega^{\otimes *})$.}
\label{Calculation of Cotor zero}

Now we begin to make cohomological calculations. The global sections (i.e., $H^0$) of the line bundles $\omega^{\otimes *}$ over $\mathcal{M}_{fmA}^{2-buds}$ are already nontrivial and interesting. The graded abelian group 
\begin{equation*}
 \Gamma(\omega^{\otimes *}, \mathcal{M}_{fmA}^{2-buds}) \cong 
 H^0_{fl}(\mathcal{M}_{fmA}^{2-buds};\omega^{\otimes *}) 
 \cong \Cotor^{0,2*}_{L^A_{2-buds}B}(L^A_{2-buds},L^A_{2-buds})\end{equation*} 
is simply the cotensor product $L^A_{2-buds}\Box_{L^A_{2-buds}B} L^A_{2-buds}$, i.e., the kernel of the difference 
\begin{align*} \eta_R - \eta_L : L^A_{2-buds} &\rightarrow L^A_{2-buds}B\end{align*}
of the unit maps on the Hopf algebroid $(L^A_{2-buds},L^A_{2-buds}B)$.
Let $A$ be a commutative ring of characteristic zero. It follows from the presentation for $(L^A_{2-buds},L^A_{2-buds}B)$ given in \eqref{iso 445010} through \eqref{iso 445014} that we have a commutative square
\[\xymatrix{
 L^A_{2-buds} \ar[r]^{\eta_R - \eta_L} \ar[d]^{\cong} & L^A_{2-buds}B \ar[d]^{\cong} \\
 \Symm_A(\Dr^A) \ar[r]^{\delta} & \Symm_A(\Dr^A)[t]
}\]
where $\delta$ is the $A$-module homomorphism which is given as follows:
\begin{itemize}
\item On the zeroth symmetric power $\Symm_A^0(\Dr^A) = A$ of $\Dr^A$, $\delta$ is the zero map.
\item On the first symmetric power $\Symm_A^1(\Dr^A) = \Dr^A$ of $\Dr^A$, $\delta$ is the {\em fundamental functional $\sigma_2$}, defined in Definition \ref{def of fundamental functional},
\[ \sigma_2: \Dr^A \rightarrow A = A\{ t\} \subseteq \Symm^0_A(\Dr^A)[t] \subseteq \Symm_A(\Dr^A)[t].\]
That is, $\delta(\gamma) = 2t$ and $\delta(c_a) = (a^2-a)t$.
\item On the second symmetric power $\Symm_A^2(\Dr^A)$ of $\Dr^A$, $\delta$ is given by
\begin{align}
\label{eq 09984} \delta(xy) &= \left((\sigma_2 x)y + x(\sigma_2 y)\right)t + (\sigma_2 x)(\sigma_2 y) t^2,
\end{align}
where $x,y\in \Dr^A$.
\item On the third symmetric power $\Symm_A^3(\Dr^A)$ of $\Dr^A$, $\delta$ is given by
\begin{align}
\label{eq 09985} \delta(xyz) &= \left((\sigma_2 x)yz + x(\sigma_2 y)z + xy(\sigma_2 z)\right)t \\
\nonumber &\ \ \ \ + \left((\sigma_2 x)(\sigma_2 y)z + (\sigma_2 x)y(\sigma_2 z) + x(\sigma_2 y)(\sigma_2 z)\right)t^2 \\
\nonumber &\ \ \ \ + (\sigma_2 x)(\sigma_2 y)(\sigma_2 z) t^3,
\end{align}
where $x,y,z\in \Dr^A$.
\item More generally\footnote{Of course the formula \eqref{eq 09986} subsumes formulas \eqref{eq 09984} and \eqref{eq 09985}. 
We have chosen to present the ideas in this redundant way, because it is easier to read formulas \eqref{eq 09984} and \eqref{eq 09985}, and to grasp the simple pattern they fit into, than to parse the general formula \eqref{eq 09986}.}, on the $n$th symmetric power $\Symm_A^n(\Dr^A)$ of $\Dr^A$, $\delta$ is given by
\begin{align}
\label{eq 09986} \delta(x_1\dots x_n) &= \sum_{i=1}^n \sum_{U\in P_i(\{ 1, \dots ,n\})} \sigma^U_2(x_1,\dots ,x_n) t^{i},\end{align}
where $P_i(\{ 1, \dots ,n\})$ is the set of $i$-element subsets of $\{ 1, \dots ,n\}$, and where $\sigma^U_2(x_1,\dots ,x_n)$ is the product of the elements $x_1 \dots x_n$ with $\sigma_2$ applied to each $x_j$ such that $j\in U$. 
\end{itemize}

Hence our task is to calculate the kernel of $\delta$. We accomplish this in Theorem \ref{cotor0 thm}. First, we need to define a certain function $\tau_n$. Suppose that $I$ is an ideal in a commutative ring $A$. Write $\iota: I\rightarrow A$ for the inclusion map. For each positive integer $n$, we can consider the $n$th symmetric power $\Symm_A^n(I)$. Given an element $x_1\cdot x_2 \cdot \dots \cdot x_n$ of $\Symm_A^n(I)$, we could take one of the elements $x_j\in I$ and regard it as an element of $A$, with the idea that $x_1 \cdot \dots x_{j-1} \cdot \iota(x_j)x_{j+1} \cdot \dots x_n$ is then an element of $\Symm_A^{n-1}(I)$. That construction would not quite yield a map $\Symm_A^n(I) \rightarrow \Symm_A^{n-1}$: the trouble is that we cannot single out a particular factor $x_j$ in the product $x_1\cdot x_2 \cdot \dots \cdot x_n$, since $x_1\cdot x_2 \cdot \dots \cdot x_n$ is an element of the {\em symmetric} power $\Symm_A^n(I)$. To get a well-defined, natural map $\Symm_A^n(I) \rightarrow \Symm_A^{n-1}(I)$, we must sum over the values $j=1, 2, \dots ,n$, yielding the $A$-module morphism
\begin{align} \label{map a0995} 
 \tau_n^I: \Symm_A^{n}(I) &\rightarrow \Symm_A^{n-1}(I) \\
\nonumber x_1\cdot\dots\cdot x_n 
  &\mapsto \iota(x_1)x_2\cdot\dots\cdot x_n \\
\nonumber  &\ \ \ \ + x_1\cdot\iota(x_2)x_3\cdot\dots\cdot x_n \\
\nonumber  &\ \ \ \ + \dots + x_1\cdot x_2\cdot\dots \cdot \iota(x_{n-1})x_n.
\end{align}
\begin{theorem}\label{cotor0 thm}
Let $A$ be a torsion-free commutative ring. 
Write $I_2$ for the universal $\mathbb{F}_2$-point-detecting ideal $(2, a^2 - a\mbox{\ for\ all\ } a\in A)\subseteq A$ in $A$, as in \cref{Point-detecting ideals...}.
Then we have an isomorphism of $A$-modules:
\begin{align*}
 H^0_{fl}(\mathcal{M}_{fmA}^{2-buds}; \omega^{\otimes n})
  &\cong \left\{ \begin{array}{ll} 
   0 &\mbox{\ if\ } n<0 \\
   A &\mbox{\ if\ } n=0 \\
   0 &\mbox{\ if\ } n=1 \\
   \ker \tau_n^{I_2} 
    &\mbox{\ if\ } n>1 ,\end{array}\right. 
\end{align*} 
where $\tau_n^{I_2}$ is the $A$-module morphism defined in \eqref{map a0995}.
\end{theorem}
\begin{proof}
By Proposition \ref{point-detecting iso}, $\sigma_2$ yields an isomorphism of $A$-modules $\tilde{\sigma_2}: \Dr^A \stackrel{\cong}{\longrightarrow} I_2$. 
The projection of $\delta(x_1 \dots x_n)$ to the $(n-1)$st symmetric power $\Symm_A^{n-1}(\Dr^A)\{ t\} \subseteq \Symm_A(\Dr^A)[t]$ is given by the formula $\delta(x_1 \dots x_n) = \tau_n^{I_2}(x_1 \dots x_n) t$. 

Consequently, for an element of $\Symm_A^n(\Dr^A)$ to be in the kernel of $\delta$, that element must be in the kernel of $\tau_n^{I_2}$. We claim that the converse is also true. That is, we claim that the kernel of $\delta\mid_{\Symm_A^n(\Dr^A)}$ is precisely the kernel of the map $\tau_n^{I_2}$. To prove this claim, suppose that $x\in \Symm_A^n(\Dr^A)$ is in the kernel of $\tau_n^{I_2}$. 
We will have cause to consider a more general class of $A$-module homomorphisms, defined as follows for each positive integer $n$ and each integer $i$ such that $0\leq i\leq n$:
\begin{align*} \bigotriangle{\!}^n_i: \Symm_A^n(\Dr^A) &\rightarrow \Symm_A^{n-i}(\Dr^A) \\
 x_1 \dots x_n &\mapsto \sum_{U\in P_{i}(\{ 1, \dots ,n\})} \sigma^U_2(x_1, \dots ,x_n) .\end{align*}
We introduce the functions $\bigotriangle{\!}^n_i$ because of the following three observations:
\begin{enumerate}
\item The formula \eqref{eq 09986} for $\delta$ is equivalent to 
\begin{align}
 \label{eq 09987} \delta(x_1 \dots x_n) &= \sum_{i=1}^n \bigotriangle{\!}^n_i(x_1 \dots x_n) t^i.
\end{align}
\item By a simple combinatorial argument, the composite 
\begin{align*} \bigotriangle{\!}^{i-j}_k\circ \bigotriangle{\!}^i_j : \Symm_A^{i}(\Dr^A) &\rightarrow \Symm_A^{i-j-k}(\Dr^A) \end{align*}
is equal to the binomial coefficient $\binom{j+k}{j}$ times the map $\bigotriangle{\!}^i_{j+k}$.
\item $\bigotriangle{\!}^n_1 = \tau_n^{I_2}$.
\end{enumerate}
Consequently $\bigotriangle{\!}^n_{i}$ is equal to a product of nonzero binomial coefficients times the composite $\bigotriangle{\!}^{i+1}_1\circ \dots \circ\bigotriangle{\!}^{n-1}_1\circ\bigotriangle{\!}^n_1$. Since $A$ is torsion-free and since $x$ was assumed to be in the kernel of $\tau_n^{I_2} = \bigotriangle{\!}^n_1$, we now have that $x$ is in the kernel of $\bigotriangle{\!}^n_i$ for each $i$. Formula \eqref{eq 09987} now yields that $x$ is in the kernel of $\delta$. Hence $\Cotor^{0,2n}_{L^A_{2-buds}B}(L^A_{2-buds},L^A_{2-buds}) \cong \ker \delta$ coincides with the kernel of $\tau_n^{I_2}$, as claimed.
\end{proof}

In particular, if $A$ is torsion-free, then $H^0_{fl}(\mathcal{M}_{fmA}^{2-buds}; \omega^{\otimes 2})$ is isomorphic to the kernel of the canonical multiplication map 
\begin{equation}\label{mult map 1} 
\Symm_A^2(I_2)\rightarrow I_2^2,
\end{equation} 
i.e., the kernel of the canonical comparison map from the symmetric square of $I_2$ to the Rees module $\Rees_A^2(I_2)$. If $A$ is a domain, then $I_2^2\subseteq A$ is torsion-free, so the kernel of the multiplication map \eqref{mult map 1} must be torsion. Torsion in $\Symm_A(I)$ is well-studied in commutative algebra: see \cite{MR563225} for a nice entry-point into the literature. The kernel of the multiplication map $\Symm_A^2(I) \rightarrow I^2$, in particular, coincides with the {\em delta-invariant} of a finitely-generated ideal $I$: see \cite{MR282964}, or Corollary 1.2 of \cite{MR610474}. Since the delta-invariant $\delta(I)$ is known to agree with the second Andre-Quillen homology group $H_2(A, A/I; A/I)$ of $A/I$ regarded as an $A$-algebra, with coefficients in $A/I$, we have:
\begin{corollary}\label{H0 cor}
Let $A$ be a Noetherian integral domain of characteristic zero. Then
the following $A$-modules are isomorphic:
\begin{itemize}
\item The sections $H^0_{fl}(\mathcal{M}_{fmA}^{2-buds}; \omega^{\otimes 2})$ of the second tensor power $\omega^{\otimes 2}$ of the bundle of invariant differentials on $\mathcal{M}_{fmA}^{2-buds}$.
\item The Andre-Quillen homology group $H_2(A, A/I_2; A/I_2)$.
\end{itemize}
\end{corollary}
Since at least the 1980s, the kernel of the map $\Symm_A^n(I) \rightarrow I^n$ has drawn attention in commutative algebra\footnote{For example, in the paper \cite{MR3563012}, it is remarked that ``Finding the defining equations of Rees rings is a classical problem in elimination theory that amounts to determining the kernel $\mathcal{A}$ of the natural map from the symmetric algebra $\Symm(I)$ onto $\mathcal{R}$.''}, especially in the case $n=2$. For example, an ideal $I$ in an integral domain $R$ is called {\em syzygetic} if the map $\Symm_A^2(I) \rightarrow I^2$ is injective (equivalently, an isomorphism); see \cite{MR690710} for some discussion and relevant results. We have:
\begin{corollary}\label{H0 cor 2}
Let $A$ be a Noetherian integral domain of characteristic zero. Then $H^0_{fl}(\mathcal{M}_{fmA}^{2-buds}; \omega^{\otimes 2})$ vanishes if and only if the universal $\mathbb{F}_2$-point-detecting ideal of $A$ is syzygetic.
\end{corollary}

If $I_2$ can be generated by a regular sequence, or more generally a $d$-sequence in the sense of Huneke \cite{MR563225}, then $\Symm_A^n(I_2) \rightarrow I_2^n$ is an isomorphism for {\em all} $n$, and consequently $H^0_{fl}(\mathcal{M}_{fmA}^{2-buds}; \omega^{\otimes n})$ vanishes for $n>0$. This explains why, in Theorem 3.2 of \cite{MR745362}, Ravenel obtained the vanishing of $\Cotor^{0,n}_{L^AB}(L^A,L^A)$ for all $n\neq 0$ and all number rings $A$: the same would be true for any regular integral domain of characteristic zero. More generally:
\begin{theorem}\label{H0 vanishing for CM}
If $A$ is a Cohen-Macaulay integral domain of characteristic zero, then $H^0_{fl}(\mathcal{M}_{fmA}^{2-buds}; \omega^{\otimes n})$ is trivial for all $n\neq 0$.
\end{theorem}
\begin{proof}
Since $A$ is Cohen-Macaulay, the ideal $I_2$ can be generated by a regular sequence, so the multiplication map
\begin{equation}\label{mult map 2}\Symm_A^n(I_2)\rightarrow I_2^n\end{equation} is injective for all $n$. The map \eqref{mult map 2} is equal to the function $\bigotriangle{\!}^n_n$ defined in the proof of Theorem \ref{cotor0 thm}. In that proof, it was shown that $\bigotriangle{\!}^n_n$ factors as a product of nonzero binomial coefficients times the composite $\bigotriangle{\!}^{2}_1\circ\dots \circ\bigotriangle{\!}^{n-1}_1\circ \bigotriangle{\!}^n_1$. Hence the injectivity of \eqref{mult map 2} implies the injectivity of the composite \begin{equation}\label{composite 093f} \bigotriangle{\!}^{i+1}_1\circ\dots \circ\bigotriangle{\!}^{n-1}_1\circ \bigotriangle{\!}^n_1\end{equation}
for each $i=1, \dots, n-1$. The composite \eqref{composite 093f} is equal to a product of nonzero binomial coefficients times $\bigotriangle{\!}^{n}_i$, so since $A$ is an integral domain of characteristic zero, $\bigotriangle{\!}^{n}_i$ is injective for all $j$. Hence, by equation \eqref{eq 09987}, the cobar complex differential morphism $\delta: \Symm_A(\Dr^A) \rightarrow\Symm_A(\Dr^A)[t]$ is injective in positive internal degrees $>0$, so $\Cotor^{0,n}_{L^A_{2-buds}B}(L^A_{2-buds},L^A_{2-buds})$ vanishes for $n>0$.
\end{proof}

\subsection{Consequences for torsion in $L^A$.}

To date, there is no known example of an integral domain $A$ such that the classifying rings $L^A$ or $L^A_{2-buds}$ have nontrivial $A$-torsion. 
The $H^0$ calculations from \cref{Calculation of Cotor zero}, together with Theorem \ref{lazard ring localization iso}, yield some insight about when and why $L^A_{2-buds}$ is torsion-free. In this section, we obtain the first known example of an integral domain $A$ such that $L^A_{2-buds}$ is {\em not} torsion-free. 

We begin with a simple observation for the case where $A$ is a field:
\begin{prop}\label{H0 vanishing for fields A}
Let $K$ be a field of characteristic zero. Then $H^0_{fl}(\mathcal{M}_{fmK}^{2-buds}; \omega^{\otimes n})$ is trivial for $n\neq 0$.
\end{prop}
\begin{proof}
Since $K$ is a field, the symmetric powers $\Sym^n_K(I_2)$ are all free $K$-modules, hence torsion-free. Consequently $\tau_n^{I_2}$ is injective. Now Theorem \ref{cotor0 thm} yields the result.
\end{proof}

\begin{prop}\label{torsion in LA2-buds}
Let $A$ be an integral domain of characteristic zero. Then, for each $n > 0$, the $A$-module $H^0_{fl}(\mathcal{M}_{fmA}^{2-buds}; \omega^{\otimes n})$ is $A$-torsion, and an $A$-submodule of $L^A_{2-buds}$.
\end{prop}
\begin{proof}
We have isomorphisms  
\begin{align*} 
 H^0_{fl}(\mathcal{M}_{fmA}^{2-buds}; \omega^{\otimes *}) 
  &\cong \Cotor^{0,2*}_{L^A_{2-buds}B}(L^A_{2-buds},L^A_{2-buds}) \\
  &\cong L^A_{2-buds}\Box_{L^A_{2-buds}B} L^A_{2-buds} \subseteq L^A_{2-buds},\end{align*}
so $H^0_{fl}(\mathcal{M}_{fmA}^{2-buds}; \omega^{\otimes n})$ is an $A$-submodule of the degree $n$ summand of $L^A_{2-buds}$. 

Write $K$ for the field of fractions of $A$. From the localization theorem (Theorem \ref{lazard ring localization iso}) and Proposition \ref{H0 vanishing for fields A}, we have isomorphisms
\begin{align*}
 0 &\cong H^0_{fl}(\mathcal{M}_{fmK}^{2-buds}; \omega^{\otimes n}) \\
   &\cong H^0_{fl}(\mathcal{M}_{fmA}^{2-buds}; \omega^{\otimes n})\otimes_A K\end{align*}
for $n\neq 0$. Hence $H^0_{fl}(\mathcal{M}_{fmA}^{2-buds}; \omega^{\otimes n})$ must be torsion for $n\neq 0$.
\end{proof}

Corollary \ref{H0 cor 2} and Proposition \ref{torsion in LA2-buds} now yield:
\begin{corollary}
Let $A$ be an integral domain of characteristic zero. If the ideal $I_2$ in $A$ is not syzygetic, then $L^A_{2-buds}$ is not torsion-free.
\end{corollary}

\begin{example}
Let $A = \mathbb{Z}[a,b,c,x]/(2a-(x^2-x)b,\ 2b - (x^2-x)c)$. Write $\hat{x}$ as shorthand for $x^2-x$. Then $\hat{x}a$ and $\hat{x}b$ and $\hat{x}c$ are elements of the ideal $I_2$ of $A$. It is routine to check that the product
\[ \hat{x}a\cdot \hat{x}c - \hat{x}b\cdot \hat{x}b \in \Sym^2_A(I_2)\]
is nonzero, and furthermore that
\begin{align*}
 2(\hat{x}a\cdot \hat{x}c - \hat{x}b\cdot \hat{x}b)
  &= \hat{x}b\cdot 2\hat{x}b - 2\hat{x}b\cdot \hat{x}b = 0.
\end{align*}
Consequently there is a nonzero $2$-torsion element of $\Sym^2_A(I_2) \cong H^0_{fl}(\mathcal{M}_{fmA}^{2-buds}; \omega^{\otimes 2})$, hence also a nonzero $2$-torsion element in $L^A_{2-buds}$. To the author's knowledge, this is the first known example of an integral domain $A$ such that $L^A_{2-buds}$ has nontrivial torsion.
\end{example}

\section{Cohomology of $\mathcal{M}_{fm\mathbb{Z}}^{2-buds}$ in low degrees.}
\label{Cohomology calculations.}

We consider the moduli stack $\mathcal{M}_{fmA}^{2-buds}$ of formal $A$-module $2$-buds in the case $A = \mathbb{Z}$. We will calculate the first cohomology group of the moduli stack of formal group $2$-buds, i.e., $\Cotor^{1,*}_{L^{\mathbb{Z}}_{2-buds}B}(L^{\mathbb{Z}}_{2-buds},L^{\mathbb{Z}}_{2-buds})$. The author does not know where this specific calculation appears in the literature, but the author does not believe that this particular calculation should be seen as especially new: it is quite similar to $2$-primary calculations of portions of the Adams-Novikov spectral sequence $E_2$-term, and also to various standard $2$-primary calculations in Iwasawa theory. Nevertheless it is worth the effort to present the calculation here, as the base case $A=\mathbb{Z}$ will be used in \cref{H0 and H1 of the moduli of...} as input for the extension-of-formal-multiplications spectral sequence which converges to $H^*_{fl}(\mathcal{M}_{fmA};\omega^{\otimes *})$ for rings $A$ other than $\mathbb{Z}$.

\subsection{Construction of the $2$-adic spectral sequence.}
\label{Construction of the 2-adic...}

The $\mathbb{Z}$-module $\Dr^{\mathbb{Z}}$ is free on the generator $\gamma$, and consequently the $\mathbb{Z}[t]$-comodule algebra $L^{\mathbb{Z}}_{2-buds}$ is isomorphic to $\mathbb{Z}[\gamma]$, with the coaction
\begin{align}
\nonumber \psi: \mathbb{Z}[\gamma] &\rightarrow \mathbb{Z}[t]\otimes_{\mathbb{Z}}\mathbb{Z}[\gamma] \\
\label{coaction 35085}  \psi(\gamma) &= 1\otimes \gamma + t\otimes 2,
\end{align}
as described above in \eqref{coaction map 1}.
The ideal $(2)$ of $L^{\mathbb{Z}}_{2-buds}$ is closed under the $\mathbb{Z}[t]$-coaction, i.e., the $(2)$-adic filtration
\[ L^{\mathbb{Z}}_{2-buds} \supseteq (2) \supseteq (2)^2 \supseteq (2)^3 \supseteq \dots \]
of $L^{\mathbb{Z}}_{2-buds}$ is a multiplicative filtration and also a filtration by subcomodules. Writing $E_0^{(2)}L^{\mathbb{Z}}_{2-buds}$ for the associated graded comodule of the $2$-adic filtration on $L^{\mathbb{Z}}_{2-buds}$, we get a multiplicative conditionally convergent spectral sequence\footnote{It is automatic that this spectral sequence converges to $\Cotor^{p,q}_{\mathbb{Z}[t]}\left(\mathbb{Z},(L_{2-buds}^{\mathbb{Z}})^{\hat{}}_2\right)$. The fact that it also converges to $\Cotor^{p,q}_{\mathbb{Z}[t]}(\mathbb{Z},L_{2-buds}^{\mathbb{Z}})^{\hat{}}_2$ is a consequence of the finite generation of $L_{2-buds}^{\mathbb{Z}}$ and of $L_{2-buds}^{\mathbb{Z}}B$ in each degree. Similar spectral sequence convergence results for rings $A$ other than $\mathbb{Z}$, and for the full moduli stack of formal $A$-modules and not merely of formal $A$-module $2$-buds, and for various filtrations including (but not limited to) the $2$-adic filtration, follow from the finiteness and completion results in \cref{finiteness properties subsection}.}
\begin{align}
\label{adic ss 1} E_1^{p,q,u} \cong \Cotor^{p,q,u}_{\mathbb{Z}[t]}(\mathbb{Z},E_0^{(2)}L_{2-buds}^{\mathbb{Z}}) 
  &\Rightarrow \Cotor^{p,q}_{\mathbb{Z}[t]}(\mathbb{Z},L_{2-buds}^{\mathbb{Z}})^{\hat{}}_2 \\
  &\cong H^p_{fl}(\mathcal{M}_{fm\mathbb{Z}};\omega^{\otimes q/2})^{\hat{}}_2 \\
\nonumber d_r : E_r^{p,q,u} &\rightarrow E_r^{p+1,q,u+r}.
\end{align} 
We will refer to spectral sequence \eqref{adic ss 1} as the {\em $2$-adic spectral sequence.} 
To be clear about the notation: the abutment $\Cotor^{*,*}_{\mathbb{Z}[t]}(\mathbb{Z},L_{2-buds}^{\mathbb{Z}})^{\hat{}}_2$ is the $2$-adic completion of $\Cotor^{*,*}_{\mathbb{Z}[t]}(\mathbb{Z},L_{2-buds}^{\mathbb{Z}})$. The spectral sequence is trivial for odd $q$, so the tensor power $\omega^{\otimes q/2}$ is well-defined.

The associated graded $\mathbb{Z}[t]$-comodule algebra $E_0^{(2)}L_{2-buds}^{\mathbb{Z}}$ of the $2$-adic filtration on $L_{2-buds}^{\mathbb{Z}}$ is isomorphic to $\mathbb{F}_2[\tilde{2},\gamma]$, with $\tilde{2}$ in $2$-adic filtration degree $1$, and with trivial coaction.
Consequently the $E_1$-term of the $2$-adic spectral sequence is isomorphic to 
\begin{align}\nonumber \Cotor^{*,*,*}_{\mathbb{Z}[t]}\left(\mathbb{Z},\mathbb{F}_2[\tilde{2},\gamma]\right) 
  &\cong \Cotor^{*,*,*}_{\mathbb{F}_2[t]}\left(\mathbb{F}_2,\mathbb{F}_2\right)\otimes_{\mathbb{F}_2}\mathbb{F}_2[\tilde{2},\gamma] \\
\label{iso f37k}  &\cong \mathbb{F}_2[\eta, P\eta, P^2\eta, P^3\eta, \dots ]\otimes_{\mathbb{F}_2}\mathbb{F}_2[\tilde{2},\gamma],\end{align}
with tridegrees $\tilde{2} \in E_1^{0,0,1},$ and $\gamma \in E_1^{0,2,0},$ and $P^j\eta \in E_1^{1,2^{j+1},0}$.
Isomorphism \eqref{iso f37k} is due to the isomorphism
\begin{align}
 \Cotor_{\mathbb{F}_2[t]}^{*,*}\left(\mathbb{F}_2,\mathbb{F}_2\right) 
\label{iso f40k}  &\cong \mathbb{F}_2[\eta, P\eta ,P^2\eta, P^3\eta, \dots ],
\end{align}
where \begin{itemize}
\item $P^n\eta\in \Cotor_{\mathbb{F}_2[t]}^{1,2^{n+1}}(\mathbb{F}_2,\mathbb{F}_2)$,
\item $P^n\eta$ is represented in the cobar complex of $(\mathbb{F}_2,\mathbb{F}_2[t])$ by the $1$-cocycle $t^{2^n}$,
\item and $P$ is the algebraic Steenrod operation $P^0$ which operates in $\Cotor$ by applying the Frobenius operation to cocycle representatives in the cobar complex. For these ideas, see the material on algebraic Steenrod operations in Appendix 1 of \cite{MR860042}, or \cite{MR0281196}. 
\end{itemize}
A straightforward way to see isomorphism \eqref{iso f40k} is to observe that $\mathbb{F}_2[t]$ splits, {\em as a coalgebra}, as the tensor product of the tensor factors $\mathbb{F}_2[t^{2^n}]/(t^{2^n})^2$ over all $n\geq 0$. The $\Cotor$-algebra of $\mathbb{F}_2[\epsilon]/\epsilon^2$ is polynomial on a single generator in cohomological degree $1$, represented by the $1$-cocycle $\epsilon$ in the cobar complex of $\mathbb{F}_2[\epsilon]/\epsilon^2$. The cobar complex multiplication is concatenation of tensors, and in the case of $\mathbb{F}_2[t]$, it is routine to make the cocycle-level calculation to verify that the coalgebra splitting 
\begin{align}  
\label{iso f41k} \mathbb{F}_2[t] 
  &\cong \bigotimes_{n\geq 0}\mathbb{F}_2[t^{2^n}]/(t^{2^n})^2 \end{align}
induces a ring isomorphism
\begin{align*}  
 \Cotor^*_{\mathbb{F}_2[t]}(\mathbb{F}_2,\mathbb{F}_2) 
  &\cong \bigotimes_{n\geq 0}\Cotor^*_{\mathbb{F}_2[t^{2^n}]/(t^{2^n})^2}(\mathbb{F}_2,\mathbb{F}_2) \end{align*}
despite \eqref{iso f41k} not respecting the ring structure.

\subsection{Running the $2$-adic spectral sequence: $H^n_{fl}(\mathcal{M}_{fm\mathbb{Z}}^{2-buds}; \omega^{\otimes *})$ for $n=0,1$.}
\label{cotor for Z}

As the $2$-adic spectral sequence is a multiplicative spectral sequence, to calculate the $d_1$-differential it suffices to calculate the $d_1$ differential on the generators $\tilde{2},\gamma,$ and $P^j\eta$ for each $j=0,1,2, \dots$.
We accomplish this using cocycle representatives for each such generator in the cobar complex for $\mathbb{Z}[t]$ with coefficients in $E_0^{(2)}L_{2-buds}^{\mathbb{Z}}$:
\begin{itemize}
\item 
Since $\tilde{2}$ is represented by the $0$-cocycle $\tilde{2}$, which lifts to the $0$-cocycle $2$ in the cobar complex for $\mathbb{Z}[t]$ with coefficients in $L_{2-buds}^{\mathbb{Z}}$, we have 
\begin{align*} d_1(\tilde{2}) &= 0,\end{align*} and in fact $\tilde{2}$ is an infinite cocycle.
\item
The generator $\gamma$ is represented by the $0$-cocycle $\gamma$, which lifts to the $0$-cochain $\gamma$ in the cobar complex for $\mathbb{Z}[t]$ with coefficients in $L_{2-buds}^{\mathbb{Z}}$. We have
\begin{align*} \delta(\gamma) &= t\otimes 2\in \mathbb{Z}[t]\otimes_{\mathbb{Z}} L_{2-buds}^{\mathbb{Z}}\end{align*}
in that cobar complex, and since $t\otimes 2$ represents a cocycle representative for $\tilde{2}\eta$ in the cobar complex for $\mathbb{Z}[t]$ with coefficients in $E_0^{(2)}L_{2-buds}^{\mathbb{Z}}$, we have 
\begin{align*} d_1(\gamma) &= \tilde{2}\eta.\end{align*}
\item 
The generator $P^j\eta$ is represented by the $1$-cocycle $t^{2^j}\otimes 1$, which lifts to the $1$-cochain $t^{2^j}\otimes 1$ in the cobar complex for $\mathbb{Z}[t]$ with coefficients in $L_{2-buds}^{\mathbb{Z}}$. We have
\begin{align*} \delta(t^{2^j}\otimes 1) 
  &= \sum_{i=1}^{2^j} \binom{2^j}{i} t^{2^j-i}\otimes t^i\otimes 1\end{align*}
in that cobar complex, whose unique term of least $2$-adic filtration is \linebreak $\binom{2^j}{2^{j-1}} t^{2^{j-1}}\otimes t^{2^{j-1}}\otimes 1$, since the central binomial coefficient $\binom{2^j}{2^{j-1}}$ has $2$-adic valuation $1$. Consequently 
\begin{align*} 
 d_1(P^j\eta) &= \frac{\binom{2^j}{2^{j-1}}}{2} \tilde{2} (P^{j-1}\eta)^2.\end{align*}
\item Consequently, by the Leibniz rule, we have 
\begin{align*}
 d_1\left( \tilde{2}^h \gamma^i P^j\eta\right)
  &= \tilde{2}^{h+1}\gamma^{i-1}\left( i\eta P^j\eta + \gamma \frac{\binom{2^j}{2^{j-1}}}{2} (P^{j-1}\eta)^2\right),
\end{align*}
with the understanding that negative powers of $\gamma$ and of $P$ are zero. Consequently the $\Cotor^1$-line in the $E_2$-page of the $2$-adic spectral sequence consists of the $\mathbb{F}_2[\tilde{2}]$-linear combinations of the elements $\gamma^{2i}\eta$ and the elements $\gamma^{2i}(\gamma\eta - P\eta)$, for $i\geq 0$. We write $Q\eta$ as shorthand for $\gamma\eta - P\eta$, i.e., the cohomology class of the $1$-cocycle $t\otimes \gamma - t^2 \otimes 1$ in the cobar complex $C_{\mathbb{Z}[t]}^{\bullet}(E_0^{(2)}L^{\mathbb{Z}}_{2-buds})$ for $\mathbb{Z}[t]$ with coefficients in $E_0^{(2)}L^{\mathbb{Z}}_{2-buds}$.
\end{itemize}
In the remaining calculations in the $2$-adic spectral sequence, we will be sloppy about the distinction between $t\otimes \gamma - t^2 \otimes 1$ and $t\otimes \gamma + t^2 \otimes 1$, since they represent the same class in the associated graded of the $2$-adic filtration.

We can calculate the $\Cotor^0$ and $\Cotor^1$ lines on later pages by similar arguments, together with the calculations
\begin{align*}
 d_2(\gamma^2) 
  &= \left[ (t\otimes 2 + 1\otimes \gamma)^2 - 1\otimes \gamma^2\right] \\
  &= \left[ 4(t^2\otimes 1 + t\otimes \gamma)\right] \\
  &\sim \tilde{2}^2Q\eta\mbox{\ in\ the\ cobar\ complex\ } C^{\bullet}_{\mathbb{Z}[t]}(E_0^{(2)}L_{2-buds}^{\mathbb{Z}}),\\
 d_3(\gamma^4) 
  &= \left[ (t\otimes 2 + 1\otimes \gamma)^4 - 1\otimes \gamma^4\right] \\
  &= \left[ 2^3(2t^4\otimes 1 + 2^2t^3\otimes \gamma + 3t^2\otimes \gamma^2 + t\otimes \gamma^3)\right] \\
  &\sim \tilde{2}^3\gamma^2Q\eta\mbox{\ in\ the\ cobar\ complex\ } C^{\bullet}_{\mathbb{Z}[t]}(E_0^{(2)}L_{2-buds}^{\mathbb{Z}}),\\
 d_4(\gamma^8)
  &\sim \tilde{2}^4\gamma^6Q\eta\mbox{\ in\ the\ cobar\ complex\ } C^{\bullet}_{\mathbb{Z}[t]}(E_0^{(2)}L_{2-buds}^{\mathbb{Z}}),
\end{align*}
and in general,
\begin{align}
\label{2adic diff formula} d_{r+1}(\gamma^{2^r}) 
  &\sim \tilde{2}^{r+1}\gamma^{2^{r}-2}Q\eta\mbox{\ in\ the\ cobar\ complex\ } C^{\bullet}_{\mathbb{Z}[t]}(E_0^{(2)}L_{2-buds}^{\mathbb{Z}}).
\end{align}
We are using the symbol $\sim$ to denote the equivalence relation ``is cohomologous to, modulo terms of higher $2$-adic filtration.''

Differential formula \eqref{2adic diff formula} yields the following description of the $\Cotor^0$ and $\Cotor^1$ lines on each page of the $2$-adic spectral sequence:
\begin{equation}
\begin{array}{l|l|l}
r    & \Cotor^0\mbox{-line in the } E_r\mbox{-page} & \Cotor^1\mbox{-line in the } E_r\mbox{-page} \\ 
\hline
1    &  \mathbb{F}_2[\tilde{2},\gamma]                    & \mathbb{F}_2[\tilde{2},\gamma]\left\{ \eta,P\eta,P^2\eta, \dots\right\}   \\
2    &  \mathbb{F}_2[\tilde{2},\gamma^2]                    & \mathbb{F}_2[\tilde{2},\gamma^2]\left\{ \eta,Q\eta\right\}/\tilde{2}\eta   \\
3    &  \mathbb{F}_2[\tilde{2},\gamma^4]                    & \frac{\mathbb{F}_2[\tilde{2},\gamma^4]\left\{ \eta,Q\eta,\gamma^2\eta,\gamma^2Q\eta\right\}}{\left(\tilde{2}\eta,\tilde{2}^2Q\eta,\tilde{2}\gamma^2\eta\right)}   \\
4    &  \mathbb{F}_2[\tilde{2},\gamma^8]                    & \frac{\mathbb{F}_2[\tilde{2},\gamma^8]\{ \eta,Q\eta,\gamma^2\eta,\gamma^2Q\eta,\gamma^4\eta,\gamma^4Q\eta,\gamma^6\eta,\gamma^6Q\eta\}}{\left(\tilde{2}\eta,\tilde{2}^2Q\eta,\tilde{2}\gamma^2\eta,\tilde{2}^3\gamma^2Q\eta,\tilde{2}\gamma^4\eta,\tilde{2}^2\gamma^4Q\eta,\tilde{2}\gamma^6\eta\right)}  ,
    \end{array}\end{equation}
and in the limit, the $\Cotor^0$-line in the $E_{\infty}$-page is $\mathbb{F}_2[\tilde{2}]$, while the $\Cotor^1$-line in the $E_{\infty}$-page is
\begin{equation*}
 \mathbb{F}_2[\tilde{2}]\left\{ \gamma^{2n}\eta,\gamma^{2n}Q\eta \ \forall n\geq 0\right\}/\left(\tilde{2}\gamma^{2n}\eta,\tilde{2}^{1+\nu_2(2n+2)}\gamma^{2n}Q\eta\ \forall n\geq 0\right).
\end{equation*}
Resolving the extension problems to pass from the $E_{\infty}$-page to the abutment $\Cotor^{*,*}_{\mathbb{Z}[t]}(\mathbb{Z},L^{\mathbb{Z}}_{2-buds})^{\hat{}}_2$, we have
\begin{align*}
 \Cotor^{0}_{\mathbb{Z}[t]}(\mathbb{Z},L^{\mathbb{Z}}_{2-buds})^{\hat{}}_2 & \cong \hat{\mathbb{Z}}_2,\\
 \Cotor^{1}_{\mathbb{Z}[t]}(\mathbb{Z},L^{\mathbb{Z}}_{2-buds})^{\hat{}}_2 & \cong \frac{\mathbb{Z}\left\{ \gamma^{2n}\eta,\gamma^{2n}Q\eta \ \forall n\geq 0\right\}}{\left(2\gamma^{2n}\eta,2^{1+\nu_2(2n+2)}\gamma^{2n}Q\eta\ \forall n\geq 0\right)}.
\end{align*}
From the $2$-adic completions of these $\Cotor$-groups, it is easy to make the cocycle-level calculations and to use the finiteness results from \cref{finiteness properties subsection} to deduce that {\em before} $2$-adic completion, we must have
\begin{align*}
 \Cotor^{0}_{\mathbb{Z}[t]}(\mathbb{Z},L^{\mathbb{Z}}_{2-buds}) & \cong \mathbb{Z},\\
 \Cotor^{1}_{\mathbb{Z}[t]}(\mathbb{Z},L^{\mathbb{Z}}_{2-buds}) & \cong \frac{\mathbb{Z}\left\{ \gamma^{2n}\eta,\gamma^{2n}Q\eta \ \forall n\geq 0\right\}}{\left(2\gamma^{2n}\eta,2^{1+\nu_2(2n+2)}\gamma^{2n}Q\eta\ \forall n\geq 0\right)}.
\end{align*}

\subsection{Running the $2$-adic spectral sequence: $H^{2}_{fl}(\mathcal{M}_{fm\mathbb{Z}}^{2-buds}; \omega^{\otimes n})$ for $n\leq 3$.}
\label{cotor2 for Z}
For later spectral sequence calculations in \cref{COER ss}, it will be useful to have calculated $\Cotor^{2,2n}_{\mathbb{Z}[t]}(\mathbb{Z},L^{\mathbb{Z}}_{2-buds})\cong H^2_{fl}(\mathcal{M}_{fm\mathbb{Z}}^{2-buds}; \omega^{\otimes n})$ for a few small values of $n$. We do this by running the $2$-adic spectral sequence in internal\footnote{The spectral sequence of a filtered cochain complex of abelian groups is bigraded: it has the cohomological degree, and the filtration degree. If the cochain complex is furthermore a filtered cochain complex of {\em graded} abelian groups, then the spectral sequence has a third grading, traditionally called the {\em internal} grading. For example, since $(L^A,L^AB)$ is a {\em graded} Hopf algebroid, the cobar complex of $(L^A,L^AB)$ is a cochain complex of {\em graded} $A$-modules, and the $2$-adic spectral sequence has an {\em internal} grading as a consequence. The internal degree is the degree $q$ in \eqref{adic ss 1}.} degrees $\leq 6$. In principle, there is no reason that similar calculations could not be done for a much wider range of internal and cohomological degrees. However, in order to keep this paper at a manageable length, we confine our attention to only the most immediately relevant calculations.

\subsubsection{Internal degrees $q < 4$.}
In these internal degrees, the $2$-adic spectral sequence has no summands which contribute to $\Cotor^2$.
\subsubsection{Internal degrees $q = 4,5$.}

In internal degree $4$, the $2$-adic spectral sequence $E_1$-page is straightforwardly calculated. The charts are as follows:
\begin{figure}[H]
\centering
\begin{minipage}{.5\textwidth}
\centering
\begin{tikzpicture}[trim left=0cm,xscale=1.2,yscale=0.7]
\draw (-0.35,3.5) -- (-0.35,-0.35) -- (2.5,-0.35);
\draw (-1,0) node{$u=0$};
\draw (-1,1) node{$u=1$};
\draw (-1,2) node{$u=2$};
\draw (-1,3) node{$u=3$};
\draw (0,-0.6) node{$p=0$};
\draw (1,-0.6) node{$p=1$};
\draw (2,-0.6) node{$p=2$};
\draw (0,0) node{$\gamma^2$};
\draw (0.8,0) node{$Q\eta$};
\draw (1.2,0) node{$\gamma\eta$};
\draw (2,0) node{$\eta^2$};
\draw (0,1) node{$\bullet$};
\draw (0,2) node{$\bullet$};
\draw (0,3) node{$\bullet$};
\draw (0.8,1) node{$\bullet$};
\draw (0.8,2) node{$\bullet$};
\draw (0.8,3) node{$\bullet$};
\draw (1.2,1) node{$\bullet$};
\draw (1.2,2) node{$\bullet$};
\draw (1.2,3) node{$\bullet$};
\draw (2,1) node{$\bullet$};
\draw (2,2) node{$\bullet$};
\draw (2,3) node{$\bullet$};
\draw[->,color=black,shorten >=0.0cm,shorten <=0.2cm] (0,0) -- (0,4) node{};
\draw[->,color=black,shorten >=0.0cm,shorten <=0.2cm] (0.8,0) -- (0.8,4) node{};
\draw[->,color=black,shorten >=0.0cm,shorten <=0.2cm] (1.2,0) -- (1.2,4) node{};
\draw[->,color=black,shorten >=0.0cm,shorten <=0.2cm] (2,0) -- (2,4) node{};
\draw[->,color=red,shorten >=0.1cm,shorten <=0.2cm] (0,0) -- (0.8,2) node{};
\draw[->,color=red,shorten >=0.1cm,shorten <=0.2cm] (1.2,0) -- (2,1) node{};
\end{tikzpicture} 
\caption{$2$-adic SS\\ \centerline{$E_1$-page,\ \ $q=4$}\\ \centerline{$d_1(\gamma\eta) = \tilde{2}\eta^2$} \\ \centerline{$d_2(\gamma^2) = \tilde{2}^2Q\eta$}}
\label{ss fig 343 2-adic q=4 E1}
\end{minipage}%
\begin{minipage}{.5\textwidth}
\centering
\begin{tikzpicture}[trim left=0cm,xscale=1.2,yscale=0.7]
\draw (-0.35,3.5) -- (-0.35,-0.35) -- (2.5,-0.35);
\draw (-1,0) node{$u=0$};
\draw (-1,1) node{$u=1$};
\draw (-1,2) node{$u=2$};
\draw (-1,3) node{$u=3$};
\draw (0,-0.6) node{$p=0$};
\draw (1,-0.6) node{$p=1$};
\draw (2,-0.6) node{$p=2$};
\draw (1.0,0) node{$Q\eta$};
\draw (2,0) node{$\eta^2$};
\draw (1.0,1) node{$\bullet$};
\draw[-,color=black,shorten >=0.0cm,shorten <=0.2cm] (1.0,0) -- (1.0,1) node{};
\end{tikzpicture} 
\caption{$2$-adic SS\\ \centerline{$E_3\cong E_{\infty}$-page,\ \ $q=4$}}
\label{ss fig 343 2-adic q=4 E3}
\end{minipage}%
\end{figure}

Empty bidegrees are understood to be trivial. Each nontrivial element name is understood to be an $\mathbb{F}_2$-linear basis element.
Vertical black arrows indicate towers of multiplications by $\tilde{2}$, so for example, the $E_1$-page is a free $\mathbb{F}_2[\tilde{2}]$-algebra on the four elements $\gamma^2, Q\eta, \gamma\eta$, and $\eta^2$. The $E_{\infty}$-page is $\mathbb{F}_2[\tilde{2}]\{ Q\eta, \eta^2\}/\left( \tilde{2}^2Q\eta, \tilde{2}\eta^2\right)$. Resolving the extensions, we have that 
\begin{align}
\nonumber \Cotor^{n,4}_{(L^{\mathbb{Z}}_{2-buds},L^{\mathbb{Z}}_{2-buds}B)}\left(L^{\mathbb{Z}}_{2-buds},L^{\mathbb{Z}}_{2-buds}\right) 
  & \cong H^n(\mathcal{M}_{fm\mathbb{Z}}^{2-buds}; \omega^{\otimes 2}) \\
\label{q=4 cotor}  &\cong \left\{ \begin{array}{ll} 
   0 &\mbox{\ if\ } n=0 \\
   \mathbb{Z}/4\mathbb{Z}\{ Q\eta\} &\mbox{\ if\ } n=1 \\
   \mathbb{Z}/2\mathbb{Z}\{ \eta^2\} &\mbox{\ if\ } n=2 \\
   0 &\mbox{\ if\ } n>2.\end{array}\right.
\end{align}

Since $L^{\mathbb{Z}}_{2-buds}$ is concentrated in even internal degrees, the $\Cotor$-groups \linebreak $\Cotor_{\mathbb{Z}[t]}^{*,*}(\mathbb{Z},L^{\mathbb{Z}}_{2-buds})$ vanish in internal degree $5$. 

\subsubsection{Internal degrees $q = 6$.}

The $2$-adic spectral sequence charts in internal degree $6$ are as follows:
\begin{figure}[H]
\centering
\begin{minipage}{.5\textwidth}
\centering
\begin{tikzpicture}[trim left=0cm,xscale=1.8,yscale=0.7]
\draw (-0.35,3.5) -- (-0.35,-0.35) -- (3.5,-0.35);
\draw (-1,0) node{$u=0$};
\draw (-1,1) node{$u=1$};
\draw (-1,2) node{$u=2$};
\draw (-1,3) node{$u=3$};
\draw (0,-0.6) node{$p=0$};
\draw (1,-0.6) node{$p=1$};
\draw (2,-0.6) node{$p=2$};
\draw (3,-0.6) node{$p=3$};
\draw (0,0) node{$\gamma^3$};
\draw (0,1) node{$\bullet$};
\draw (0,2) node{$\bullet$};
\draw (0,3) node{$\bullet$};
\draw (0.8,0) node{$\gamma^2\eta$};
\draw (0.8,1) node{$\bullet$};
\draw (0.8,2) node{$\bullet$};
\draw (0.8,3) node{$\bullet$};
\draw (1.2,0) node{$\gamma Q\eta$};
\draw (1.2,1) node{$\bullet$};
\draw (1.2,2) node{$\bullet$};
\draw (1.2,3) node{$\bullet$};
\draw (1.8,0) node{$\eta Q\eta$};
\draw (1.8,1) node{$\bullet$};
\draw (1.8,2) node{$\bullet$};
\draw (1.8,3) node{$\bullet$};
\draw (2.2,0) node{$\gamma\eta^2$};
\draw (2.2,1) node{$\bullet$};
\draw (2.2,2) node{$\bullet$};
\draw (2.2,3) node{$\bullet$};
\draw (3,0) node{$\eta^3$};
\draw (3,1) node{$\bullet$};
\draw (3,2) node{$\bullet$};
\draw (3,3) node{$\bullet$};
\draw[->,color=black,shorten >=0.0cm,shorten <=0.2cm] (0,0) -- (0,4) node{};
\draw[->,color=black,shorten >=0.0cm,shorten <=0.2cm] (0.8,0) -- (0.8,4) node{};
\draw[->,color=black,shorten >=0.0cm,shorten <=0.2cm] (1.2,0) -- (1.2,4) node{};
\draw[->,color=black,shorten >=0.0cm,shorten <=0.2cm] (1.8,0) -- (1.8,4) node{};
\draw[->,color=black,shorten >=0.0cm,shorten <=0.2cm] (2.2,0) -- (2.2,4) node{};
\draw[->,color=black,shorten >=0.0cm,shorten <=0.2cm] (3,0) -- (3,4) node{};
\draw[->,color=red,shorten >=0.1cm,shorten <=0.2cm] (0,0) -- (0.8,1) node{};
\draw[->,color=red,shorten >=0.1cm,shorten <=0.2cm] (1.2,0) -- (1.8,1) node{};
\draw[->,color=red,shorten >=0.1cm,shorten <=0.2cm] (2.2,0) -- (3,1) node{};
\end{tikzpicture} 
\caption{$2$-adic SS\\ \centerline{$E_1$-page,\ \ $q=6$} \\ \centerline{$d_1(\gamma^3) = \tilde{2} \gamma^2\eta$} \\ \centerline{$d_1(Q\eta \cdot \gamma) = \tilde{2}\eta\cdot Q\eta$} \\ \centerline{$d_1(\eta^2 \cdot \gamma) = \tilde{2}\eta^3$}}
\label{ss fig 343 2-adic q=6 E1}
\end{minipage}%
\end{figure}
\begin{figure}[H]
\begin{minipage}{.5\textwidth}
\centering
\begin{tikzpicture}[trim left=0cm,xscale=1.8,yscale=0.7]
\draw (-0.35,3.5) -- (-0.35,-0.35) -- (3.5,-0.35);
\draw (-1,0) node{$u=0$};
\draw (-1,1) node{$u=1$};
\draw (-1,2) node{$u=2$};
\draw (-1,3) node{$u=3$};
\draw (0,-0.6) node{$p=0$};
\draw (1,-0.6) node{$p=1$};
\draw (2,-0.6) node{$p=2$};
\draw (3,-0.6) node{$p=3$};
\draw (1.0,0) node{$\gamma^2\eta$};
\draw (2.0,0) node{$\eta Q\eta$};
\draw (3,0) node{$\eta^3$};
\end{tikzpicture} 
\caption{$2$-adic SS \\ $E_3\cong E_{\infty}$-page,\ \ $q=6$}
\label{ss fig 343 2-adic q=6 E2}
\end{minipage}%
\end{figure}

There is no room for nontrivial extensions, so from the $E_{\infty}$-page we have that 
\begin{align*}
 \Cotor^{n,6}_{(L^{\mathbb{Z}}_{2-buds},L^{\mathbb{Z}}_{2-buds}B)}\left(L^{\mathbb{Z}}_{2-buds},L^{\mathbb{Z}}_{2-buds}\right) 
  & \cong H^n(\mathcal{M}_{fm\mathbb{Z}}^{2-buds}; \omega^{\otimes 3}) \\
  &\cong \left\{ \begin{array}{ll} 
   0 &\mbox{\ if\ } n=0 \\
   \mathbb{Z}/2\mathbb{Z}\{ \gamma^2 \eta\} &\mbox{\ if\ } n=1 \\
   \mathbb{Z}/2\mathbb{Z}\{ \eta\cdot Q\eta\} &\mbox{\ if\ } n=2 \\
   \mathbb{Z}/2\mathbb{Z}\{ \eta^3 \} &\mbox{\ if\ } n=3 \\
   0 &\mbox{\ if\ } n>3.\end{array}\right.
\end{align*}

\begin{remark}\label{anss remark}
As a consequence of the calculations in this section, we have that the bigraded ring $\coprod_{s,t} H^s_{fl}(\mathcal{M}_{fm\mathbb{Z}}^{2-buds}; \omega^{\otimes t})$ is $2$-locally isomorphic in the range $t-s\leq 3$ to $\coprod_{s,t} H^s_{fl}(\mathcal{M}_{fm\mathbb{Z}}; \omega^{\otimes t})$, i.e., the input for the Adams-Novikov spectral sequence. Drawn with the Adams conventions:

\begin{minipage}{1.0\textwidth}
\centering
\begin{tikzpicture}[trim left=0cm,xscale=1.8,yscale=0.7]
\clip (-1.5,-1.5) rectangle (6.5,6.5);
\fill[green!50!white] (6.8,0.0) -- (0.0,6.8) -- (6.8,6.8) -- cycle;
\draw (-0.35,6.5) -- (-0.35,-0.35) -- (6.5,-0.35);
\draw (-1,0) node{$s=0$};
\draw (-1,1) node{$s=1$};
\draw (-1,2) node{$s=2$};
\draw (-1,3) node{$s=3$};
\draw (-1,4) node{$s=4$};
\draw (-1,5) node{$s=5$};
\draw (-1,6) node{$s=6$};
\draw (0,-0.6) node{$t-s=0$};
\draw (1,-0.6) node{$t-s=1$};
\draw (2,-0.6) node{$t-s=2$};
\draw (3,-0.6) node{$t-s=3$};
\draw (4,-0.6) node{$t-s=4$};
\draw (5,-0.6) node{$t-s=5$};
\draw (6,-0.6) node{$t-s=6$};
\draw (0,0) node{$\mathbb{Z}\{ 1\}$};
\draw (1,1) node{$\mathbb{Z}/2\mathbb{Z}\{ \eta\}$};
\draw (2,2) node{$\mathbb{Z}/2\mathbb{Z}\{ \eta^2\}$};
\draw (3,3) node{$\mathbb{Z}/2\mathbb{Z}\{ \eta^3\}$};
\draw (3,1) node{$\mathbb{Z}/4\mathbb{Z}\{ Q\eta\}$};
\draw (4,2) node{$\mathbb{Z}/2\mathbb{Z}\{ \eta\cdot Q\eta\}$};
\draw (5,1) node{$\mathbb{Z}/2\mathbb{Z}\{ \gamma^2\eta\}$};
\end{tikzpicture} 
\end{minipage}

The green-shaded region lies beyond what we have just calculated. 
Of course it is possible to run the $2$-adic spectral sequence in higher internal degrees for a more far-reaching comparison of $H^*_{fl}(\mathcal{M}_{fm\mathbb{Z}}^{2-buds}; \omega^{\otimes *})$ with the $E_2$-term $H^*_{fl}(\mathcal{M}_{fm\mathbb{Z}}; \omega^{\otimes *})$ of the Adams-Novikov spectral sequence, but in this paper our priority is on making calculations of $H^*_{fl}(\mathcal{M}_{fm\mathbb{A}}^{2-buds}; \omega^{\otimes *})$ in low bidegrees for a very wide range of rings $A$. Those calculations begin in the next section, and the low-degree calculations of $H^*_{fl}(\mathcal{M}_{fm\mathbb{Z}}^{2-buds}; \omega^{\otimes *})$ here are only in the service of those in the next section.
\end{remark}

\section{$H^0$ and $H^1$ of the moduli of formal $A$-module $2$-buds.}
\label{H0 and H1 of the moduli of...}

\subsection{The symmetric filtration on symmetric powers.}

In \cref{COER ss}, we will construct and use a spectral sequence which passes from the cohomology of the moduli stack of formal $\mathbb{Z}$-module $2$-buds to the cohomology of the moduli stack of formal $A$-module $2$-buds. The construction of this spectral sequence relies on a certain filtration of symmetric powers. This subsection is about that filtration. The filtration is constructed in Definition \ref{def of symm filt}, and its associated graded is calculated in Proposition \ref{assoc gr of symm filt}.

The filtration is, at least under certain hypotheses, quite well-known. Given a commutative ring $R$ and a short exact sequence of $R$-modules
\begin{equation}\label{ses 430005} 0 \rightarrow M^{\prime} \rightarrow M^{} \rightarrow M^{\prime\prime} \rightarrow 0,\end{equation}
we ask for an increasing $R$-module filtration on the symmetric algebra $\Symm_R^*(M)$ such that the associated graded $R$-algebra $E^0\Symm_R^*(M)$ is isomorphic to \linebreak $\Symm_R^*(M^{\prime})\otimes_R \Symm_R^*(M^{\prime\prime})$. It is standard that such a filtration exists, {\em if the $R$-modules $M,M^{\prime},M^{\prime\prime}$ are projective}; see exercise 5.16 in chapter II of \cite{MR0463157}, for example.

However, to build the spectral sequence that we will use in \cref{COER ss}, it will be necessary to relax those hypotheses slightly. The essential condition is that each of the maps $s_1, s_2, \dots $ in a certain sequence, \eqref{symm layer sequence}, are one-to-one. It turns out that this condition is satisfied in the case of interest in \cref{COER ss}, even though certain of the $R$-modules involved are not projective, and even though the short exact sequence \eqref{ses 430005} will not split.

The results in this section are elementary, but technical, involving colimits over certain ``truncated-cube-shaped'' diagrams. The author apologizes for not being able to find a simpler way to present the ideas. Surely these ideas cannot be new, and must be well-known within some circles, but we were unable to find a reference in the literature.

Now we begin the relevant definitions. First we must introduce the indexing categories for certain colimits.
\begin{definition}\label{def of cube cats}\leavevmode
\begin{itemize}
\item Let $\mathcal{I}$ denote the category with
\begin{itemize}
\item two objects, $0$ and $1$,
\item a single homomorphism $0\rightarrow 1$,
\item and no non-identity endomorphisms.
\end{itemize}
In other words, $\mathcal{I}$ is the partially-ordered set $\{0,1\}$, regarded as a category.
\item Let $n$ be a nonnegative integer.
Let $\mathcal{I}^n$ be the $n$-fold Cartesian product category $\mathcal{I}\times \dots \times\mathcal{I}$. That is, $\mathcal{I}^n$ is the partially-ordered set of $n$-tuples $\{0,1\} \times \dots \times \{ 0,1\}$, regarded as a category. The relevant partial ordering is the one in which $(a_1, \dots ,a_n)\leq (b_1, \dots ,b_n)$ if and only if $a_i\leq b_i$ for all $i$.
\item Let $i,n$ be nonnegative integers, with $i\leq n$. Let $\mathcal{I}^n_i$ denote the full subcategory of $\mathcal{I}^n$ containing precisely those objects $(a_1, \dots ,a_n)$ such that $\sum_j a_j \leq i$. 
\end{itemize}
\end{definition}
For example, when $n=2$, we have the following pictures of $\mathcal{I}^2_0,\mathcal{I}^2_1,$ and $\mathcal{I}^2_2 = \mathcal{I}^2$:
\begin{figure}[H]
\centering
\begin{minipage}{.30\textwidth}
\centering
\[ \xymatrix{ (0,0) & \\ & 
}\]
\caption{$\mathcal{I}^2_0$}
\end{minipage}%
\begin{minipage}{.30\textwidth}
\centering
\[ \xymatrix{ (0,0)\ar[r]\ar[d] & (1,0) \\ (0,1) & 
}\]
\caption{$\mathcal{I}^2_1$}
\end{minipage}%
\begin{minipage}{.40\textwidth}
\centering
\[ \xymatrix{ (0,0)\ar[r]\ar[d] & (1,0) \ar[d] \\ (0,1)\ar[r] & (1,1) 
}\]
\caption{$\mathcal{I}^2_2 = \mathcal{I}^2$}
\end{minipage}%
\end{figure}

When $n=3$, we have the following pictures:
\begin{figure}[H]
\centering
\begin{minipage}{.5\textwidth}
\centering
\[ \xymatrix{ (0,0,0) & \\ & 
}\]
\caption{$\mathcal{I}^3_0$}
\end{minipage}%
\begin{minipage}{.5\textwidth}
\centering
\[ \xymatrix{ (0,0,0)\ar[r]\ar[d]\ar[rdd] & (1,0,0)  \\ (0,1,0) &  \\ & (0,0,1) 
}\]
\caption{$\mathcal{I}^3_1$}
\end{minipage}%
\end{figure}
\begin{figure}[H]
\begin{minipage}{.5\textwidth}
\centering
\[ \xymatrix{ (0,0,0)\ar[r]\ar[d]\ar[rdd] & (1,0,0)\ar[d]\ar[rdd] & \\ (0,1,0)\ar[r]\ar[rdd] & (1,1,0) & \\ & (0,0,1)\ar[r]\ar[d] & (1,0,1) \\ & (0,1,1) &
}\]
\caption{$\mathcal{I}^3_2$}
\end{minipage}%
\begin{minipage}{.5\textwidth}
\centering
\[ \xymatrix{ (0,0,0)\ar[r]\ar[d]\ar[rdd] & (1,0,0)\ar[d]\ar[rdd] & \\ (0,1,0)\ar[r]\ar[rdd] & (1,1,0)\ar[rdd] & \\ & (0,0,1)\ar[r]\ar[d] & (1,0,1)\ar[d] \\ & (0,1,1) \ar[r] & (1,1,1)
}\]
\caption{$\mathcal{I}^3_3 = \mathcal{I}^3$}
\end{minipage}%
\end{figure}

\begin{definition}\label{def of res and L}
Let $R$ be a commutative ring, let $n$ be a nonnegative integer, let $M_0,M_1$ be $R$-modules, and let $f: M_0 \rightarrow M_1$ be an $R$-module homomorphism.
\begin{itemize}
\item
Let $\mathcal{F}: \mathcal{I}^n \rightarrow \Mod(R)$ be the functor given by sending $(a_1, \dots ,a_n)$ to $M_{a_1}\otimes_R \dots \otimes_R M_{a_n}$.
\item
For each nonnegative integer $i\leq n$, we have the restriction $\mathcal{F}\mid_{\mathcal{I}^n_i}: \mathcal{I}^n_i \rightarrow \Mod(R)$ of $\mathcal{F}$ to the full subcategory $\mathcal{I}^n_i$ of $\mathcal{I}^n$. 
\item 
For each nonnegative integer $i < n$, the restriction functor $\res_i: \Mod(R)^{\mathcal{I}^n_{i+1}} \rightarrow \Mod(R)^{\mathcal{I}^n_{i}}$ admits a left adjoint. We write $L_i$ for this left adjoint.
\end{itemize}
\end{definition}
Of course $\colim \mathcal{F}$ is simply the $n$-fold tensor power of $M_1$, and using the natural action of $\Sigma_n$ on $\mathcal{I}_n$, we have $(\colim\mathcal{F})_{\Sigma_n} \cong \Symm^n_R(M_1)$. The idea of introducing the subcategories $\mathcal{I}^n_i$ of $\mathcal{I}^n$ is to obtain a useful filtration of the symmetric power $\Symm^n_R(M_1)$.

\begin{lemma}\label{symm filt quotients lemma}
Let $R,n,M_0,M_1,f,\mathcal{F}$ be as in Definition \ref{def of res and L}. 
Let $\tilde{M}_0$ denote $M_0$, and let $\tilde{M}_1$ denote the cokernel of $f: M_0 \rightarrow M_1$. 
For each nonnegative integer $i<n$, the cokernel of the counit map
\begin{equation}\label{counit 0349}  L_i\res_i(\mathcal{F}\mid_{\mathcal{I}^n_{i+1}}) \rightarrow \mathcal{F}\mid_{\mathcal{I}^n_{i+1}} \end{equation}
is the functor $\mathcal{I}^n_{i+1} \rightarrow \Mod(R)$ that sends $(a_1, \dots ,a_n)$ to $0$ if $\sum_j a_j \leq i$, and sends $(a_1, \dots ,a_n)$ to $\tilde{M}_{a_1} \otimes \dots \otimes \tilde{M}_{a_n}$ if $\sum_j a_j = i+1$.
\end{lemma}
\begin{proof}
By the pointwise formula for Kan extensions (classical; see \cite{MR1712872} for example), $L_i\res_i(\mathcal{F}\mid_{\mathcal{I}^n_{i+1}})$ is given as follows:
\begin{align}\label{counit 0350} L_i\res_i(\mathcal{F}\mid_{\mathcal{I}^n_{i+1}})(a_1, \dots ,a_n) &= \left\{\begin{array}{ll} 
  \mathcal{F}(a_1, \dots ,a_n) &\mbox{\ if\ } \sum_j a_j \leq i \\
  \colim_{(b_1, \dots ,b_n)<(a_1, \dots ,a_n)} \mathcal{F}(b_1, \dots ,b_n) &\mbox{\ if\ } \sum_j a_j = i+1.\end{array}\right.
\end{align}
using the partial ordering on $\mathcal{I}^n$ from Definition \ref{def of cube cats}. 

Cokernels in functor categories are computed levelwise, so the fact that \eqref{counit 0350} coincides with $\mathcal{F}$ if $\sum_ja_j\leq i$ tells us that the cokernel $c_i: \mathcal{I}^n_{i+1} \rightarrow \Mod(R)$ of the map \eqref{counit 0349} vanishes on all tuples $(a_1, \dots ,a_n)$ such that $\sum_j a_j \leq i$, as claimed.

As for those tuples $(a_1, \dots ,a_n)$ such that $\sum_j a_j = i+1$: let $\mathcal{I}_{<(a_1, \dots ,a_{n})}$ denote the full subcategory of $\mathcal{I}^n_{i+1}$ consisting of those tuples $(b_1, \dots ,b_{n})$ satisfying $(b_1, \dots ,b_{n})<(a_1, \dots ,a_{n})$.
Then $\mathcal{I}_{<(a_1, \dots ,a_{n})}$ is isomorphic to $\mathcal{I}^{i+1}_{i}$. Let $\tilde{\mathcal{F}}_{(a_1, \dots, a_n)}: \mathcal{I}_{<(a_1, \dots ,a_{n})} \rightarrow \Mod(R)$ denote the constant functor taking the value $\mathcal{F}(a_1, \dots,a_n)$. We have a natural map $\mathcal{F}\mid_{\mathcal{I}_{<(a_1, \dots ,a_{n})}}\rightarrow \tilde{\mathcal{F}}_{(a_1, \dots, a_n)}$ which is an isomorphism when evaluated on $(a_1,\dots ,a_n)$. Again using the fact that cokernels are calculated levelwise in functor categories, the value of the cokernel of \eqref{counit 0349} at $(a_1, \dots, a_n)$ agrees with the value of the cokernel of the composite map
\begin{equation}\label{counit 0351}  L_i\res_i(\mathcal{F}\mid_{\mathcal{I}^n_{i+1}})\mid_{\mathcal{I}_{<(a_1, \dots ,a_{n})}}\rightarrow \mathcal{F}\mid_{\mathcal{I}_{<(a_1, \dots ,a_{n})}}\rightarrow \tilde{\mathcal{F}}_{(a_1, \dots, a_n)}
\end{equation}
at $(a_1, \dots ,a_n)$.

The cokernel of the composite map \eqref{counit 0351} is the functor which sends $(b_1, \dots ,b_n)$ to the cokernel of the map
$\mathcal{F}(b_1, \dots ,b_n) \rightarrow\mathcal{F}(a_1, \dots,a_n)$. In the case that $\sum_j b_j = i$, this cokernel is precisely $\tilde{M}_{b_1} \otimes_R \dots \otimes_R \tilde{M}_{b_n}$. It is routine to verify that the colimit of the cokernel of \eqref{counit 0351} is consequently $\tilde{M}_{a_1} \otimes_R \dots \otimes_R \tilde{M}_{a_n}$. Consequently the cokernel of \eqref{counit 0351} sends $(a_1, \dots ,a_n)$ to $\tilde{M}_{a_1} \otimes_R \dots \otimes_R \tilde{M}_{a_n}$, as claimed.\end{proof}

\begin{definition}\label{def of symm filt}
Let $R,n,M_0,M_1,f,\mathcal{F}$ be as in Definition \ref{def of res and L}. 
Let $n$ be a nonnegative integer.
\begin{itemize}
\item Given a nonnegative integer $i<n$, write $\tilde{\res}_i$ for the restriction functor $\Mod(R)^{\mathcal{I}^n} \rightarrow \Mod(R)^{\mathcal{I}^n_i}$, i.e., the composite of the functors $\res_{i}\circ \dots \circ \res_{n-2}\circ\res_{n-1}$ from Definition \ref{def of res and L}. Write $\tilde{L}_i$ for its left adjoint, i.e., the composite of the functors $\res_{n-1}\circ \dots \circ \res_{i+1}\circ L_i$, also from Definition \ref{def of res and L}.
\item Given a positive integer $i<n$, we have the natural transformation of functors $\mathcal{I}^n \rightarrow \Mod(R)$
\begin{equation}\label{counit map 0353}  \tilde{L}_{i-1}\tilde{\res}_{i-1}\mathcal{F} \rightarrow \tilde{L}_i\tilde{\res}_i\mathcal{F}. \end{equation}
We write $s_i$ for the induced map of $R$-modules
\begin{equation}\label{counit map 0354} \left(\colim \tilde{L}_{i-1}\tilde{\res}_{i-1}\mathcal{F}\right)_{\Sigma_n} \rightarrow \left(\colim \tilde{L}_i\tilde{\res}_i\mathcal{F}\right)_{\Sigma_n} .\end{equation}
By the {\em symmetric layer sequence} we mean the sequence of $R$-module maps
\begin{equation}\label{symm layer sequence} \xymatrix{ \Symm_R^n(M_0) \ar[r]^= & \left(\colim \tilde{L}_0\tilde{\res}_0\mathcal{F}\right)_{\Sigma_n} \ar[d]^{s_1} & \\ & \left(\colim \tilde{L}_1\tilde{\res}_1\mathcal{F}\right)_{\Sigma_n} \ar[d]^{s_{2}}  & \\ & \vdots\ar[d]^{s_{n-1}}  & \\ &\left(\colim \tilde{L}_{n-1}\tilde{\res}_{n-1}\mathcal{F}\right)_{\Sigma_n} \ar[d]^{s_n}  & \\ & \left(\colim\mathcal{F}\right)_{\Sigma_n} \ar[r]^= & \Symm_R^n(M_1).}\end{equation}
\item If each of the maps $s_i$ is injective, then the symmetric layer sequence \eqref{symm layer sequence} is a filtration of $\Symm_R^n(M_1)$, and we call this filtration the {\em symmetric filtration.}
\end{itemize}
\end{definition}

Using Lemma \ref{symm filt quotients lemma} to identify the cokernels of the maps in \eqref{counit map 0353}, we have:
\begin{prop}\label{assoc gr of symm filt}
Let $R,n,M_0,M_1,f,\mathcal{F}$ be as in Definition \ref{def of res and L}.  Suppose that each of the maps $s_i$ in the symmetric layer sequence \eqref{symm layer sequence} is injective. Then the symmetric filtration is an increasing filtration of $\Symm_R^n(M_1)$ whose associated graded $R$-module is isomorphic to the direct sum \[ \coprod_{i=0}^n \Symm_R^i(M_0)\otimes_R \Symm_R^{n-i}(\coker f).\]

If each of the maps $s_i$ in the symmetric layer sequence \eqref{symm layer sequence} is injective {\em for all nonnegative integers $n$,} then, in particular, $f$ is injective, so we regard $M_0$ as a submodule of $M_1$, and we have an isomorphism of graded $R$-algebras
\[ E_0\Symm_R^*(M_1) \cong \Symm_R^*(M_0)\otimes_R \Symm_R^*(M_1/M_0).\]
\end{prop}

\subsection{The extension-of-formal-multiplications (EFM) spectral sequence.}\label{COER ss}
We now use the calculations from \cref{Cohomology calculations.} to obtain calculations of \linebreak $H^1_{fl}(\mathcal{M}_{fmA}^{2-buds};\omega^{\otimes *})$ for a broad class of rings $A$, not just the base case $A = \mathbb{Z}$. The main tool is a spectral sequence which allows us to pass from $H^*_{fl}(\mathcal{M}_{fm\mathbb{Z}}^{2-buds};\omega^{\otimes *})$ to $H^*_{fl}(\mathcal{M}_{fmA}^{2-buds};\omega^{\otimes *})$. As this is a matter of passing from the cohomology of the moduli of formal group law $2$-buds with a small ring of formal multiplications (namely, $\mathbb{Z}$) to the cohomology of the moduli of formal group $2$-buds with a larger ring of formal multiplications (namely, $A$), we call this spectral sequence the {\em extension-of-formal-multiplications spectral sequence,} or for short, ``EFM spectral sequence.''

\begin{theorem}\label{change of endo-ring ss}
Let $A$ be a torsion-free commutative ring. Let $A/2$ denote the reduction of the ring $A$ modulo the ideal $(2)$.
Let $\tilde{\Omega}^A$ denote the free bigraded $A/2$-module\footnote{The $A/2$-module $\tilde{\Omega}^2$ is a module of ``twisted K\"{a}hler differential forms'' in the sense of \cite{MR1975301}. This is intriguing, but we do not know of any general theorems about modules of twisted K\"{a}hler differentials which give us any leverage here. Perhaps this is a reasonable direction for later investigations.} on the set of generators $\{ c_a : a\in A\}$ modulo the relations
\begin{align*}
 c_{a+b} &=c_a+c_b\mbox{\ \ for\ all\ } a,b\in A,\\
 c_{ab} &= ac_b + b^2c_a\mbox{\ \ for\ all\ } a,b\in A,\\
 c_a &=0\mbox{\ \ for\ all\ } a\in \mathbb{Z}.
\end{align*}
Then there exists a conditionally convergent spectral sequence\footnote{The author admits to finding it difficult to visualize the spectral sequence's $E_1$-term merely from the description given in \eqref{sseq 443442}. We find the charts drawn below, starting in Figure \ref{ss fig 34409}, much more helpful for visualizing the spectral sequence.} 
\begin{align}
\label{sseq 443442} E_1^{p,q,u} &\cong 
  \left\{ \begin{array}{ll} 
   A\otimes_{\mathbb{Z}} \Cotor^{p,q}_{\mathbb{Z}[t]}(\mathbb{Z},L_{2-buds}^{\mathbb{Z}})
    &\mbox{\ if\ } u=0 \\
  \coprod_{i\geq 0} \Cotor^{p,q-2(i+u)}_{\mathbb{F}_2[t]}\left(\mathbb{F}_2, \mathbb{F}_2\right)\otimes_{\mathbb{F}_2} \Symm_{A/2}^{u}(\tilde{\Omega}^{A}) \{ \gamma^{i}\} 
    &\mbox{\ if\ } u>0 \end{array}\right. \\
\nonumber  &\Rightarrow \Cotor^{p,q}_{A[t]}(A,L_{2-buds}^A) \\
\nonumber &\cong \Cotor^{p,q}_{L_{2-buds}^AB}(L_{2-buds}^A,L_{2-buds}^A)\\
\nonumber d_r: E_r^{p,q,u} &\rightarrow E_r^{p+1,q,u-r}
\end{align}
with tridegrees as follows: 
\begin{equation*}
\begin{array}{l|lllll}
\mbox{Coh.\ class}          & \mbox{Coh.\ degree\ } (p) & \mbox{Int.\ degree\ } (q)  & \mbox{Filt.\ degree\ } (u)\\
\hline \gamma                           & 0                         & 2                          & 0  \\
c_a                         & 0                         & 2                          & 1  \\
P^n \eta                    & 1                         & 2^{n+1}                     & 0, \end{array}\end{equation*} 
where the elements $\eta, P\eta, P^2\eta, \dots$ are the generators of $\Cotor^{1,*}_{\mathbb{F}_2[t]}(\mathbb{F}_2,\mathbb{F}_2)$, as in \eqref{iso f40k} in \cref{cotor for Z}, and $\eta$ is also the generator of $\Cotor^{1,*}_{\mathbb{Z}[t]}(\mathbb{Z},L_{2-buds}^{\mathbb{Z}})$, as in \cref{cotor for Z}.
\end{theorem}
\begin{proof}
Consider the symmetric layer sequence, as in \eqref{symm layer sequence}, arising from the $A$-module homomorphism $A\otimes_{\mathbb{Z}}\Dr^{\mathbb{Z}} \hookrightarrow \Dr^A$. Since $A\otimes_{\mathbb{Z}}\Dr^{\mathbb{Z}}$ is a free $A$-module of rank $1$, each of its symmetric powers is also a free $A$-module of rank $1$, and the natural map 
\[ \Symm_A^n(A\otimes_{\mathbb{Z}}\Dr^{\mathbb{Z}}) \rightarrow \Symm_A^n(\Dr^A)\]
is simply the inclusion of the $A$-submodule of $(L^A)^{2n}$ generated by $\gamma^n$. Hence the maps in the symmetric layer sequence are injective. The resulting symmetric filtration of $\Symm_A^*(\Dr^A) \cong L^A_{2-buds}$ is the increasing filtration 
\[ F_0L_{2-buds}^A \subseteq F_1L_{2-buds}^A \subseteq F_2L_{2-buds}^A \subseteq \dots\]
on $L_{2-buds}^A$ given by letting $F_0L_{2-buds}^A$ be the $A$-subalgebra of $L_{2-buds}^A$ generated by $\gamma$, and letting $F_nL_{2-buds}^A$ be the $F_0L_{2-buds}^A$-submodule of $L_{2-buds}^A$ generated by all products of up to $n$ elements $c_a$ with $a\in A$. By the coaction formulas given in \eqref{coaction map 1}, this is a filtration by $A[t]$-subcomodules of $L_{2-buds}^A$. The filtration is multiplicative, exhaustive, complete, and separated, hence we have a conditionally convergent multiplicative spectral sequence
\begin{align}
\label{sseq 443443} E_1^{p,q,u} \cong \Cotor^{p,q}_{A[t]}\left(A,F_uL_{2-buds}^A/F_{u-1}L_{2-buds}^A\right)
 &\Rightarrow \Cotor^{p,q}_{A[t]}(A,L_{2-buds}^A) \\
\nonumber d_r: E_r^{p,q,u} &\rightarrow E_r^{p+1,q,u-r}.
\end{align}

Since $\tilde{\Omega}^A \cong \Dr^A/(A\otimes_{\mathbb{Z}}\Dr^{\mathbb{Z}})$, by Proposition \ref{assoc gr of symm filt} the associated graded $A[t]$-comodule $E^0L_{2-buds}^A = \coprod_u F_uL_{2-buds}^A/F_{u-1}L_{2-buds}^A$ of the symmetric filtration is isomorphic to 
\begin{equation}\label{symm filt iso 120} \Symm_A^*(A\otimes_{\mathbb{Z}}\Dr^{\mathbb{Z}}) \otimes_A \Symm_A^*(\tilde{\Omega}^A). \end{equation}
To more explicitly identify the $E_1$-term of spectral sequence \eqref{sseq 443443}, we need to know what happens when we apply the functor $\Cotor_{A[t]}^*(A,-)$ to \eqref{symm filt iso 120}. In general there is no particularly nice K\"{u}nneth-like formula for $\Cotor$ applied to a tensor product of comodules, but in this case, we are fortunate: if $u>0$, then the $A[t]$-coaction on $F_uL_{2-buds}^A$ lands in $F_{u-1}L_{2-buds}^A$. Consequently the $A[t]$-coaction on $\coprod_{u>0}F_uL_{2-buds}^A/F_{u-1}L_{2-buds}^A$ is trivial.
Hence, for positive $u$, we use Proposition \ref{assoc gr of symm filt} to obtain isomorphisms
\begin{align*}
 \Cotor^{p,q}_{A[t]}\left(A,\frac{F_uL_{2-buds}^A}{F_{u-1}L_{2-buds}^A}\right) 
  &\cong \Cotor^{p,q}_{A[t]}\left(A, \coprod_{i\geq 0} \Symm_A^{i}(A\otimes_{\mathbb{Z}} \Dr^{\mathbb{Z}})\otimes_A \Symm_A^{u}(\tilde{\Omega}^{A})\right) \\
  &\cong \Cotor^{p,q}_{A[t]}\left(A, \coprod_{i\geq 0} \Symm_{A}^{i}(A\otimes_{\mathbb{Z}} \Dr^{\mathbb{Z}})\otimes_A \Symm_{A/2}^{u}(\tilde{\Omega}^{A})\right) \\
 \\
  &\cong \left( \coprod_{i\geq 0} \Cotor^{p,*}_{A/2[t]}\left(A/2, A/2\right)\right. \\   &\ \ \ \ \ \ \left.\otimes_{A/2} \left( A/2\otimes_{\mathbb{Z}} \Symm_{\mathbb{Z}}^{i}(\Dr^{\mathbb{Z}})\right)\otimes_{A/2} \Symm_{A/2}^{u}(\tilde{\Omega}^{A})\right)^q \\
  &\cong \left(\coprod_{i\geq 0} \Cotor^{p,*}_{\mathbb{F}_2[t]}\left(\mathbb{F}_2, \mathbb{F}_2\right) \otimes_{\mathbb{Z}} \Symm_{\mathbb{Z}}^{i}(\Dr^{\mathbb{Z}})\otimes_{\mathbb{Z}} \Symm_{A/2}^{u}(\tilde{\Omega}^{A}) \right)^q.
\end{align*}
\end{proof}

To avoid any potential confusion: in the statement of Theorem \ref{cotor1 calc thm}, $A/I_2^2$ denotes the quotient of $A$ by the square of its ideal $I_2$.
\begin{theorem}\label{cotor1 calc thm}
Let $A$ be a torsion-free commutative Noetherian\footnote{The assumption that $A$ is Noetherian is used only to ensure that $I_2$ is finitely generated, so that the results on the delta-invariant from \cite{MR282964} and \cite{MR992766} apply. If $A$ is not assumed Noetherian, the theorem holds as stated, except for the identification of the summand $H^0_{fl}(\mathcal{M}_{fmA}^{2-buds}; \omega^{\otimes 2})$ in terms of the delta-invariant.} ring. 
Then we have isomorphisms of $A$-modules
\begin{align*}
 H^s_{fl}(\mathcal{M}^{2-buds}_{fmA}; \mathcal{O}) 
  &\cong \left\{ \begin{array}{ll} 
   A &\mbox{\ if\ } s=0 \\ 
   0 &\mbox{\ if\ } s\neq 0, \end{array}\right. \\
 H^s_{fl}(\mathcal{M}^{2-buds}_{fmA}; \omega) 
  &\cong \left\{ \begin{array}{ll} 
   A/I_2 &\mbox{\ if\ } s=1 \\ 
   0 &\mbox{\ if\ } s\neq 1, \end{array}\right. \\
 H^s_{fl}(\mathcal{M}^{2-buds}_{fmA}; \omega^{\otimes 2}) 
  &\cong \left\{ \begin{array}{ll} 
   \delta(I_2) &\mbox{\ if\ } s=0 \\ 
   A/I_2^2 &\mbox{\ if\ } s=1 \\ 
   A/I_2  &\mbox{\ if\ } s=2 \\
   0 &\mbox{\ otherwise},\end{array}\right. 
\end{align*}
where $\delta(I_2)$ is the delta-invariant of the universal $\mathbb{F}_2$-point-detecting ideal $I_2$ of $A$, i.e., $\delta(I_2)$ is isomorphic to the Andre-Quillen homology group $H_2(A, A/I; A/I)$, as explained in Corollary \ref{H0 cor}.
\end{theorem}
The claim that $H^1_{fl}(\mathcal{M}^{2-buds}_{fmA}; \omega^{\otimes n})$ vanishes for $n\leq 0$ is quite trivial, since the Hopf algebroid $(L^A_{2-buds},L^A_{2-buds}B)$ vanishes in negative degrees, and is easily seen to have no coalgebroid primitives in degree zero. It is the calculation of $H^1_{fl}(\mathcal{M}^{2-buds}_{fmA}; \omega)$ and of $H^1_{fl}(\mathcal{M}^{2-buds}_{fmA}; \omega^{\otimes 2})$ which takes a bit more work. We carry out this calculation with by running spectral sequence \eqref{sseq 443443} in the necessary bidegrees. 

Throughout the calculation, we will often begin with an element $a\in A$, and then need to consider the element $a^2-a\in A$. We adopt the following notation, which serves to streamline the discussion: given an element $a\in A$, we will write $\hat{a}$ for the element $a^2-a\in A$. It will become convenient later to have taken note of a few basic properties of $\hat{a}$:
\begin{align}
\nonumber \widehat{ab}  &= \hat{a}\hat{b} + a^2b + ab^2, \mbox{\ \ and} \\
\nonumber \hat{a}c_b + \hat{b}c_a &= 2c_{ab} - 2bc_a - 2ac_b \\
\label{eq 99013}  &= \hat{a}\hat{b}\gamma.
\end{align}

The EFM spectral sequence differentials preserve the internal degree $q$, so it is convenient to carry out the further calculations in the spectral sequence by proceeding one internal degree at a time. 

\subsubsection{Internal degree $q=0$.}
In negative degrees and in odd internal degrees, everything is trivial, since $L^A$ and $L^AB$ vanish in those internal degrees. In internal degree $q=0$, the spectral sequence is as depicted:

\begin{figure}[H]
\centering
\begin{minipage}{.5\textwidth}
\centering
\begin{tikzpicture}[trim left=0cm,xscale=1.2,yscale=0.7]
\draw (-0.35,3.5) -- (-0.35,-0.35) -- (2.5,-0.35);
\draw (-1,0) node{$u=0$};
\draw (-1,1) node{$u=1$};
\draw (-1,2) node{$u=2$};
\draw (-1,3) node{$u=3$};
\draw (0,-0.6) node{$p=0$};
\draw (1,-0.6) node{$p=1$};
\draw (2,-0.6) node{$p=2$};
\draw (0,0) node{$A\{1\}$};
\end{tikzpicture} 
\caption{EFM SS\\ \centerline{$E_1\cong E_{\infty}$-page,\ \ $q=0$} \\ \centerline{No differentials.}}
\end{minipage}
\end{figure}

\smallskip

Bidegrees left blank are understood to be zero. 
Consequently we have \[ \Cotor_{A[t]}^{0,0}(A,L^A_{2-buds}) \cong
H^0_{fl}(\mathcal{M}_{fmA}^{2-buds}; \mathcal{O}) \cong A\{1\}\] trivially, for degree reasons.

\subsubsection{Internal degree $q=2$.}
The EFM spectral sequence is as depicted:
\begin{figure}[H]
\centering
\begin{minipage}{.5\textwidth}
\centering
\begin{tikzpicture}[trim left=0cm,xscale=1.2,yscale=0.7]
\draw (-0.35,3.5) -- (-0.35,-0.35) -- (2.5,-0.35);
\draw (-1,0) node{$u=0$};
\draw (-1,1) node{$u=1$};
\draw (-1,2) node{$u=2$};
\draw (-1,3) node{$u=3$};
\draw (0,-0.6) node{$p=0$};
\draw (1,-0.6) node{$p=1$};
\draw (2,-0.6) node{$p=2$};
\draw (0,1) node{$\tilde{\Omega}_A$};
\draw (1,0) node{$A/2\{\eta\}$};
\draw[->,color=red,shorten >=0.3cm,shorten <=0.3cm] (0,1) -- (1,0) node{};
\end{tikzpicture} 
\caption{EFM SS\\ \centerline{$E_1$-page,\ \ $q=2$} \\ \centerline{$d_1(c_a) = \hat{a}\eta$}}
\label{ss fig 34409}
\end{minipage}%
\begin{minipage}{.5\textwidth}
\centering
\begin{tikzpicture}[trim left=0cm,xscale=1.2,yscale=0.7]
\draw (-0.35,3.5) -- (-0.35,-0.35) -- (2.5,-0.35);
\draw (-1,0) node{$u=0$};
\draw (-1,1) node{$u=1$};
\draw (-1,2) node{$u=2$};
\draw (-1,3) node{$u=3$};
\draw (0,-0.6) node{$p=0$};
\draw (1,-0.6) node{$p=1$};
\draw (2,-0.6) node{$p=2$};
\draw (0,1) node{$$};
\draw (1,0) node{$A/I_2\{\eta\}$};
\end{tikzpicture} 
\caption{EFM SS\\ \centerline{$E_2\cong E_{\infty}$-page,\ \ $q=2$}}
\label{ss fig 34409a}
\end{minipage}
\end{figure}

For degree reasons, $d_r(\eta) = 0$ for all $r$. By the coaction map \eqref{coaction map 1}, we have the differential $d_1(c_a) = \hat{a}\eta$ in the spectral sequence; this is the differential drawn in red in the $q=2$ diagram in Figure \ref{ss fig 34409}.
Consequently $\Cotor_{A[t]}^{1,2}(A,L^A_{2-buds})$ is a free $A/I_2$-module on the generator $\eta$. We also have that $\Cotor_{A[t]}^{0,2}(A,L^A_{2-buds})$ is the kernel of that $d_1$-differential, but this $\Cotor$-group was calculated already in \cref{Calculation of Cotor zero}.

\subsubsection{Internal degree $q=4$.}
We use isomorphism \eqref{q=4 cotor}, obtained from the $2$-adic spectral sequence, to identify the $u=0$-line in the EFM spectral sequence $E_1$-term:
\begin{figure}[H]
\centering
\begin{minipage}{.7\textwidth}
\centering
\begin{tikzpicture}[trim left=0cm,xscale=3.0,yscale=1.0]
\draw (-0.35,3.5) -- (-0.35,-0.35) -- (2.5,-0.35);
\draw (-1,0) node{$u=0$};
\draw (-1,1) node{$u=1$};
\draw (-1,2) node{$u=2$};
\draw (-1,3) node{$u=3$};
\draw (0,-0.6) node{$p=0$};
\draw (1,-0.6) node{$p=1$};
\draw (2,-0.6) node{$p=2$};
\draw (0,2) node{$\Symm_{A/2}^2(\tilde{\Omega}_A)$};
\draw (0,1) node{$\tilde{\Omega}_A\{\gamma\}$};
\draw (1,1) node{$\tilde{\Omega}_A\{\eta\}$};
\draw (1,0) node{$A/4\{Q\eta\}$};
\draw (2,0) node{$A/2\{\eta^2\}$};
\draw[->,color=red,shorten >=0.6cm,shorten <=0.6cm] (1,1) -- (2,0) node{};
\draw[->,color=red,shorten >=0.6cm,shorten <=0.6cm] (0,1) -- (1,0) node{};
\draw[->,color=blue,shorten >=0.6cm,shorten <=0.6cm] (0,2) -- (1,0) node{};
\end{tikzpicture} 
\caption{EFM SS $E_1$-page,\ \ $q=4$ \\ \centerline{$d_1(c_a\eta) = \hat{a}\eta^2$\ \ \ \ \ \ \ \ \ $d_1(c_a\gamma) = 2\hat{a}Q\eta$} \\ \centerline{$d_1(c_ac_b) = 0$\ \ \ \ \ \ \ \ \ $d_2(c_ac_b) = \hat{a}\hat{b}Q\eta$}}
\end{minipage}
\end{figure}
\begin{figure}[H]
\centering
\begin{minipage}{.7\textwidth}
\centering
\begin{tikzpicture}[trim left=0cm,xscale=3.0,yscale=1.0]
\draw (-0.35,3.5) -- (-0.35,-0.35) -- (2.5,-0.35);
\draw (-1,0) node{$u=0$};
\draw (-1,1) node{$u=1$};
\draw (-1,2) node{$u=2$};
\draw (-1,3) node{$u=3$};
\draw (0,-0.6) node{$p=0$};
\draw (1,-0.6) node{$p=1$};
\draw (2,-0.6) node{$p=2$};
\draw (0,2) node{$\delta(I_2)$};
\draw (0,1) node{$$};
\draw (1,1) node{};
\draw (1,0) node{$A/I_2^2\{Q\eta\}$};
\draw (2,0) node{$A/I_2\{\eta^2\}$};
\end{tikzpicture} 
\caption{EFM SS $E_{\infty}$-page,\ \ $q=4$}
\end{minipage}
\end{figure}

The differential $d_1: \Symm_{A/2}^2(\tilde{\Omega}_A) \rightarrow \tilde{\Omega}_A\{\eta\}$ vanishes due to the equalities
\begin{align}
\nonumber d_1(c_ac_b)
  &= (\hat{a} c_b + \hat{b}c_a) \eta\\
\nonumber  &= ((a^2 c_b + b c_a) + (a c_b + b^2 c_a)) \eta \\
\nonumber  &= (c_{ab} + c_{ab}) \eta \\
\label{eq 491}  &= 0\end{align}
in the associated graded $E^0L_{2-buds}^A$ of the symmetric filtration on $L_{2-buds}^A$. 
Meanwhile, by the Leibniz rule, the differential $d_1: \tilde{\Omega}_A\{\eta\} \rightarrow A\otimes_{\mathbb{Z}}\Cotor^{2,4}_{\mathbb{Z}[t]}(\mathbb{Z},L_{2-buds}^{\mathbb{Z}})$ is merely $\eta$ times the $q=2$ $d_1$-differential $\tilde{\Omega}_A \rightarrow A\otimes_{\mathbb{Z}}\Cotor^{1,2}_{\mathbb{Z}[t]}(\mathbb{Z},L_{2-buds}^{\mathbb{Z}})$. Since $t\otimes t\otimes 1$ is not a coboundary in the cobar complex of $\mathbb{Z}[t]$ with coefficients in $L_{2-buds}^{\mathbb{Z}}$, $\eta^2$ is nonzero, although $2\eta^2=0$. Hence the kernel of $d_1: E_1^{1,4,1}\rightarrow E_1^{2,4,0}$ is the kernel of the $A$-module homomorphism
\begin{align*}
 \tilde{\Omega}_A\{ \eta\} & \rightarrow A/2\{ \eta^2\} \\
 c_a\eta &\mapsto \hat{a}\eta^2.
\end{align*}
Since $d_1(c_ac_b)$ was already shown to vanish, the kernel of $d_1: E_1^{1,4,1}\rightarrow E_1^{2,4,0}$ is $E_2^{1,4,1}$. Hence we have an isomorphism
\begin{align*}
 E_2^{0,2,1} &\stackrel{\cong}{\longrightarrow} E_2^{1,4,1} \\
 x &\mapsto x\eta.
\end{align*}
Since $E_2^{0,2,1} \cong \Cotor^{0,2}_{L^A_{2-buds}B}(L^A_{2-buds},L^A_{2-buds})$ is trivial by Theorem \ref{cotor0 thm}, $E_2^{1,4,1}$ also is trivial.

We have the $d_2$-differential $E_2^{0,4,2}\rightarrow E_2^{1,4,0}$ given by the cobar complex calculation
\begin{align}
\nonumber d_2: \Symm^2_{A/2}(\tilde{\Omega}_A) \rightarrow E_2^{1,4,0}  \\
  d_2(c_ac_b) 
   &= \left[ (c_a\otimes 1 + \hat{a}\otimes t)(c_b\otimes 1 + \hat{b}\otimes t) - c_ac_b\otimes 1\right] \\
\nonumber   &= \left[ (\hat{a}c_b + \hat{b}c_a)\otimes t + \hat{a}\hat{b}\otimes t^2\right] \\
\label{eq 94013}   &= \hat{a}\hat{b}\left[ \gamma\otimes t + 1 \otimes t^2\right] \\
\nonumber   &= \hat{a}\hat{b}Q\eta,
\end{align}
with \eqref{eq 94013} a consequence of \eqref{eq 99013}.
Hence we have isomorphisms
\begin{align*}
 \Cotor_{A[t]}^{1,4}(A,L_{2-buds}^A) 
  &\cong E_3^{1,4,0} \\
  &\cong E_{\infty}^{1,4,0}\\
  &\cong A/I_2^2\{ Q\eta\}.\end{align*}

As a corollary of Remark \ref{remark on local zeta} and Theorem \ref{cotor1 calc thm}, we have:
\begin{corollary}\label{local zeta-function cor}
Let $A$ be a Noetherian integral domain of characteristic zero. Then $H^1_{fl}(\mathcal{M}_{fmA}^{2-buds}; \omega)$ is a finite abelian group of order equal to $2^{N_1}$, where $N_1$ is the number of $\mathbb{F}_2$-points of $\Spec A$, i.e., the logarithmic derivative of the $2$-local zeta-function $Z(\Spec A,t)$ evaluated at $t=0$. 
\end{corollary}

\bibliography{/home/asalch/texmf/tex/salch}{}
\bibliographystyle{plain}
\end{document}